\documentclass[11pt,a4paper]{article}
\usepackage[letterpaper,top=2cm,bottom=2cm,left=3cm,right=3cm,marginparwidth=1.75cm]{geometry}
\usepackage{amsmath}
\usepackage{graphicx}
\usepackage{array}

\newcolumntype{P}[1]{>{\centering\arraybackslash}p{#1}}
\usepackage{dsfont}
\usepackage{multirow}
\usepackage{makecell}
\usepackage{epstopdf}
\usepackage{subfigure}
\usepackage{algorithm}
\usepackage{algpseudocode}
\usepackage{cite}
\usepackage{multicol}
\usepackage{multirow}
\usepackage{mathtools}
\usepackage{amssymb}
\usepackage{color}
\usepackage{bm}
\usepackage{indentfirst}
\usepackage{url}
\usepackage{nomencl}
\usepackage{booktabs}
\usepackage{cite}
\usepackage[numbers,sort&compress]{natbib}
\usepackage{graphicx}
\usepackage{subfigure}
\usepackage{algorithm}
\usepackage{algpseudocode}
\usepackage{cite}
\usepackage{multicol}
\usepackage{multirow}
\usepackage{mathtools}
\usepackage{amssymb}
\usepackage{color}
\usepackage{bm}
\usepackage{indentfirst}
\usepackage{url}
\usepackage{nomencl}
\usepackage{booktabs}
\usepackage{cite}
\usepackage[numbers,sort&compress]{natbib}
\usepackage{amsthm}
\newtheorem{Theorem}{Theorem}
\newtheorem{Lemma}{Lemma}

\newtheorem{Remark}{Remark}

\newtheorem{Assumption}{Assumption}
\allowdisplaybreaks

\begin{document}

\title{Duality-Based Distributed Optimization with Communication Delays in Multi-Cluster Networks} 


\author{Jianzheng Wang, Guoqiang~Hu 
\thanks{This work was supported in part by Singapore Economic Development Board under EIRP grant S14-1172-NRF EIRP-IHL, and in part by the Republic of Singapore's National Research Foundation under its Campus for Research Excellence and Technological Enterprise (CREATE) Programme through a grant to the Berkeley Education Alliance for Research in Singapore (BEARS) for the Singapore-Berkeley Building Efficiency and Sustainability in the Tropics (SinBerBEST) Program.}
	\thanks{Jianzheng Wang and Guoqiang Hu are with the School of Electrical and Electronic Engineering, Nanyang Technological University, Singapore, 639798 e-mail: (wang1151@e.ntu.edu.sg, gqhu@ntu.edu.sg).}
}
\maketitle

\begin{abstract}                          
In this work, we consider solving a distributed optimization problem (DOP) in a multi-agent network with multiple agent clusters. In each cluster, the agents manage separable cost functions composed of possibly non-smooth components and aim to achieve an {agreement} on a common decision of the cluster. The global cost function is considered as the sum of the individual cost functions associated with affine {coupling} constraints on the {clusters'} decisions. To solve this problem, the dual problem is formulated by the concept of Fenchel conjugate. Then an {asynchronous distributed dual proximal gradient} (Asyn-DDPG) algorithm is proposed based on a cluster-based {partial} and {mixed} consensus protocol, by which the agents are only required to communicate with their {neighbors} with {communication delays}. An ergodic convergence result is provided, and the feasibility of the proposed algorithm is verified by solving a social welfare optimization problem in the simulation.
\end{abstract}


\section{Introduction}

\subsection{Background and Motivation}

{D}{istributed} optimization problems (DOPs) have drawn much attention in the recent few years due to their wide applications in practical problems, such as task assignment in multi-robot networks and machine learning problems, where the optimal solution of the whole network can be achieved only with local communications \cite{luo2014provably,lee2017speeding}. To adapt to the arising computational burden and privacy preservation issues in the applications with large-scale data sets, DOPs with multiple agent clusters were explored recently \cite{guo2017distributed}. Generally speaking, a multi-cluster network corresponds to the conventional multi-agent network where each node is broken into a cluster of sub-nodes such that the overall computation task for the node can be separated effectively. In a multi-cluster network, the state updates and information exchanges are usually managed at two levels: cluster level and system level, which bring more complications into the algorithm development. On the other hand, asynchronous optimization problems were actively studied recently due to the communication delays in the practical systems and the increasing demand in the computational accuracy of various fields. Motivated by the discussions above, we aim to develop an asynchronous optimization algorithm for DOPs by considering the communication delays in a multi-cluster network.

\subsection{Literature Review}

In this work, we consider solving a novel DOP with a separable cost function and general affine coupling constraints in a multi-cluster network. The proposed DOP contains {two} objectives: {intra-cluster} consensus and {network-wide} optimization. In particular, the agreement of each cluster is achieved by the consensus protocol among the agents involved. Then the optimal solution to the DOP is achieved with the presence of the coupling constraints across the clusters. Therefore, the problem formulation in this work is significantly different from the cluster-based DOPs {without} coupling constraints studied in \cite{guo2017distributed,li2019gossip,shi2019multi} in terms of both {algorithm development} and {potential applications}. Meanwhile, due to the {heterogeneous} objectives within and across the clusters, it may be challenging to directly apply the existing algorithms {without} clusters \cite{sirb2018decentralized,ramaswamy2022optimization,cao2021decentralized,tian2020achieving,tsianos2012distributed,zhang2019asyspa,chang2016asynchronous,hale2017asynchronous,bedi2018asynchronous,su2022convergence,cao2020constrained,bedi2019asynchronous,wang2022composite}.
In addition, different from \cite{guo2017distributed,li2019gossip,shi2019multi}, we consider the communication delays between the neighboring agents. To this end, an {asynchronous} optimization algorithm is proposed, where each agent is only required to access the {historical} states of its neighbors.
Specifically, by augmenting the decision variables, the original cluster-based optimization problem is reconstructed as an affine-constrained problem, whose key features are {threefold}. {First}, different from the decentralized optimization algorithms proposed in \cite{sirb2018decentralized,ramaswamy2022optimization,cao2021decentralized,tian2020achieving,tsianos2012distributed,zhang2019asyspa,chang2016asynchronous}, general affine coupling constraints are considered. {Second}, different from the problems with coupling constraints in \cite{hale2017asynchronous,bedi2018asynchronous,su2022convergence,cao2020constrained}, a distributed computation manner is realized in a general incomplete communication graph. {Third}, the asynchronous models in \cite{bedi2019asynchronous,wang2022composite} assume that the variables in the same coupling constraint should be managed by {neighboring} agents, which is a special case of general coupling constraints.

To optimize the problem of interest, the dual problem is studied in this work.
Then to realize distributed computations in general incomplete graphs, a cluster-based {partial} (in the sense that only a {subset} of the overall variables is consensual) and {mixed} (with the presence of {both} cluster-wide and network-wide agreements) consensus protocol is established, which can be more general than the existing consensus schemes studied in \cite{cao2021decentralized,sirb2018decentralized,tian2020achieving,tsianos2012distributed,zhang2019asyspa,chang2016asynchronous}, where only the conventional {global consensus} is considered. To this end, a primal-dual algorithm framework is employed in this work, where {both} the primal and dual variables (in terms of the {dual} problem) belonging to different clusters can be transmitted between the {neighboring} agents with delays. Differently, to address the asynchrony, some {central} nodes are established in \cite{su2022convergence,bedi2018asynchronous,chang2016asynchronous,hale2017asynchronous,feyzmahdavian2016asynchronous,nedic2001distributed} and {complete} communication graphs are required by \cite{cao2020constrained,hannah2018unbounded,wang2023composite,zhou2018distributed}, which may be restrictive in some large-scale networks. For instance, the primal-dual optimization scheme in \cite{hale2017asynchronous} requires certain central coordinator to update the dual variables. In \cite{cao2020constrained}, both the primal and dual variables require the information of {all} the agents via a complete graph. In view of the formulated {nonsmooth} dual functions, the existing algorithms only applicable to smooth cost functions may not be directly applied in this work \cite{tian2020achieving,hale2017asynchronous,hannah2018unbounded,sirb2018decentralized,wang2021asynchronous,su2022convergence,ramaswamy2022optimization}.

The contributions of this work are summarized as follows.

\begin{enumerate}
  \item We consider a class of DOPs in a {multi-cluster} network with composite cost functions and general affine {coupling} constraints. As a new feature, the proposed problem requires a consensual decision within each cluster and also the network-wide optimization across different clusters. By augmenting the decision variable of each cluster, the dual problem is formulated based on Fenchel conjugate.
  \item To solve the formulated dual problem, an asynchronous distributed dual proximal gradient (Asyn-DDPG) algorithm is proposed, which can address the communication delays in the network. The optimal solution can be achieved when the agents make updates in a {distributed} computation manner by accessing the {outdated} information of neighbors. An ergodic convergence rate can be obtained in terms of the function value and consensus constraint violation. The feasibility of the proposed algorithm is verified by solving a social welfare optimization problem in the simulation.
\end{enumerate}


$\mathbb{N}$ denotes the non-negative integer space.
Let $\mid {A}\mid$ be the size of set ${A}$. $\mathbb{R}_+^n$ denotes the $n$-dimensional Euclidian space only with non-negative elements.
Operator $(\cdot)^{\top}$ represents the transpose of a matrix.
Define $(a)^+ \triangleq \max\{0,a\}$ with $a$ being a scalar. $A \times B$ is the Cartesian product of sets $A$ and $B$.
$\prod_{n \in B} {A}_n$ and $\bigcap_{n \in B} {A}_n$ denote the Cartesian product and interaction of sets from ${A}_1$ to ${A}_{|B|}$, respectively. $\mathcal{P}_{A}[\cdot]$ is a Euclidean projection onto set $A$,
and $\mathbf{relint}{A}$ represents the relative interior of set ${A}$. Define $\| \mathbf{a}\|^2_{\mathbf{E}} \triangleq \mathbf{a}^{\top} \mathbf{E}\mathbf{a}$ with $\mathbf{E}$ being a square matrix. $\otimes$ is Kronecker product. $\tau_{\mathrm{max}}(\mathbf{E})$ and $\tau_{\mathrm{min}}(\mathbf{E})$ denote the largest and smallest eigenvalues of matrix $\mathbf{E}$, respectively. $\| \mathbf{E}\|_2$ and $\| \mathbf{a}\|$ refer to the spectral norm of matrix $\mathbf{E}$ and $l_2$-norm of vector $\mathbf{a}$, respectively. $\mathbf{I}_n$ is an $n$-dimensional identity matrix and $\mathbf{O}_{n \times m}$ is an $(n \times m)$-dimensional zero matrix. $\mathbf{1}_n$ and $\mathbf{0}_n$ denote the $n$-dimensional column vectors with all elements being 1 and 0, respectively. $\mathbf{E} \succ 0$ ($\mathbf{E} \succeq 0$) means that $\mathbf{E}$ is positive definite (semi-definite).
Define a diagonal-like matrix $\mathrm{diag}[\mathbf{E}_1,...,\mathbf{E}_m] \triangleq \left[
                                   \begin{array}{ccc}
                                     \mathbf{E}_1 &  & \mathbf{O} \\
                                       & \ddots &  \\
                                     \mathbf{O} &   & \mathbf{E}_m \\
                                   \end{array}
                                 \right]$,
where the $(i,j)$th element of $\mathbf{E}_n$ is the $(\sum_{l=0}^{n-1} \alpha_{l} +i, \sum_{l=0}^{n-1} \beta_{l}+j)$th element of $\mathrm{diag}[\mathbf{E}_1,...,\mathbf{E}_m]$, $ \mathbf{E}_n \in \mathbb{R}^{\alpha_{n} \times \beta_{n}}$, $\alpha_{0} = \beta_{0} = 0$, $n=1,2,...,m$.
Let $p: \mathbb{R}^n \rightarrow (-\infty,+\infty]$ be a proper, convex, and closed function. The proximal mapping of $p$ is defined by $\mathbf{prox}^{\alpha}_{p} [\mathbf{m}] \triangleq \arg \min_{\mathbf{n}} ( p(\mathbf{n}) + \frac{1}{2{\alpha}} \| \mathbf{n} - \mathbf{m}\|^2 )$, ${\alpha}>0$, $\mathbf{m} \in \mathbb{R}^n$.

\section{Problem Formulation}\label{sa3}

The considered network model, problem formulation, and relevant assumptions are introduced as follows.

\subsection{Network Model}
\subsubsection{Graphs} Consider a multi-agent network ${{G}}$, which is composed of cluster set $V \triangleq \{1,...,N\}$. Cluster $i \in V$ is defined by ${G}_i \triangleq \{{V}_i,{E}_i\}$ with agent set ${V}_i \triangleq \{1,...,n_i\}$ and undirected edge set $E_i \triangleq \{e_1,e_2,...,e_{|E_i|}\}  \subseteq \{(j,l)| j \in {V}_i, l\in {V}_i, j \neq l\}$ (no self-loop). Alternatively, the index of the $j$th agent in  cluster $i$ can also be {relabeled} by $n_{ij} \triangleq \sum_{l=0}^{i-1}n_l + j$ with $n_0=0$, {i.e.}, the agents are relabeled from cluster $1$ to $N$ according to the index in each respective cluster. Then $G$ can be described by agent set $\widehat{V} \triangleq \{1,...,\sum_{l=1}^{N}n_l\}$ and undirected edge set $E \triangleq \{\widehat{e}_1,\widehat{e}_2,...,\widehat{e}_{|E|}\} \subseteq \{(j,l)| j \in \widehat{V}, l \in \widehat{V}, j \neq l\}$ (no self-loop). Let ${{V}}_i^j \triangleq  \{ l | (j,l) \in {E}_i\}$ and $\widehat{V}^k \triangleq  \{ l | (k,l) \in {E}\}$ be the neighbor sets of the $j$th agent in $G_i$ and the $k$th agent (relabeled index) in $G$, respectively. 

\subsubsection{Indexing Law of Edges} With a given indexing law of the vertices in $V_i$, the index of edges in $E_i$ can be decided as follows. For any two distinct edges $e_k = (k_1,k_2) \in E_i$ and $e_v=(v_1,v_2) \in E_i$, if $\min\{k_1,k_2\} > \min\{v_1,v_2\}$, then $k > v$, and vice versa. For the case $\min\{k_1,k_2\} = \min\{v_1,v_2\}$, if $\max\{k_1,k_2\} > \max\{v_1,v_2\}$, then $k > v$, and vice versa. Similarly, based on the {relabeled} vertex indices, the index of edges in $E$ is decided as follows. For any two distinct edges $\widehat{e}_{\widehat{k}} = (\widehat{k}_1,\widehat{k}_2) \in E$ and $\widehat{e}_{\widehat{v}}=(\widehat{v}_1,\widehat{v}_2) \in E$, if $\min\{\widehat{k}_1,\widehat{k}_2\} > \min\{\widehat{v}_1,\widehat{v}_2\}$, then $\widehat{k} > \widehat{v}$, and vice versa. For the case $\min\{\widehat{k}_1,\widehat{k}_2\} = \min\{\widehat{v}_1,\widehat{v}_2\}$, if $\max\{\widehat{k}_1,\widehat{k}_2\} > \max\{\widehat{v}_1,\widehat{v}_2\}$, then $\widehat{k} > \widehat{v}$, and vice versa.



\subsubsection{Matrices} Let $\mathbf{L}^i \in \mathbb{R}^{n_i \times n_i}$ be the Laplacian matrix of $G_i$, where the $(j,l)$th element is defined by $[\mathbf{L}^i]_{jl} \triangleq \left\{
\begin{array}{ll}
  |V^j_i| & \hbox{if $j=l \in  V_i$} \\
  -1 & \hbox{if $(j,l) \in E_i$} \\
  0 & \hbox{otherwise}
\end{array}
\right.$ \cite{chung1997spectral}.
The graph $G_i$ can also be described by an incidence matrix $\mathbf{G}_i \in \mathbb{R}^{n_i \times |E_i|}$, where
$[\mathbf{G}_i]_{jk} \triangleq \left\{
\begin{array}{ll}
  1 & \hbox{if $e_k=(j,l) \in E_i$ and $j<l$} \\
  -1 & \hbox{if $e_k=(j,l) \in E_i$ and $j>l$} \\
  0 & \hbox{otherwise}
\end{array}
\right.$ \cite{wang2023composite}.
For a singleton graph $G_i$ (one vertex and zero edge), we define $\mathbf{G}_i \in \mathbb{R}^{1 \times 0}$ for the convenience of presentation \cite{diestel2005graph}. Similarly, we define an incidence matrix $\widehat{\mathbf{G}} \in \mathbb{R}^{|\widehat{V}| \times |E|}$ for $G$, where $[\widehat{\mathbf{G}}]_{jk} \triangleq \left\{
\begin{array}{ll}
  1 & \hbox{if $\widehat{e}_k=(j,l) \in E$ and $j<l$} \\
  -1 & \hbox{if $\widehat{e}_k=(j,l)\in E$ and $j>l$} \\
  0 & \hbox{otherwise}
\end{array}
\right.$.

\subsection{The Optimization Problem}

Let $F(\mathbf{x}) \triangleq \sum_{i \in {V}} F_{i}(\mathbf{x}_{i})$ be the global cost function of network $G$ and let $F_i(\mathbf{x}_i) \triangleq \sum_{j \in {V}_i} (f_{ij}(\mathbf{x}_{i}) + g_{ij}(\mathbf{x}_{i}))$ be the cost function of cluster $i$, where $f_{ij}+g_{ij}$ is the cost function of the $j$th agent in cluster $i$, $\mathbf{x}_{i} \in \mathbb{R}^M$, $\mathbf{x} \triangleq [\mathbf{x}^{\top}_{1},...,\mathbf{x}^{\top}_{N}]^{\top} \in \mathbb{R}^{{N}M}$.
The optimization problem of network $G$ is formulated as
\begin{align}
  \mathrm{(P1)} \quad  \min\limits_{\mathbf{x}} \quad  & \sum_{i \in {V}} \sum_{j \in {V_i}} (f_{ij}(\mathbf{x}_{i}) + g_{ij} (\mathbf{x}_{i})) \nonumber \\
   \hbox{subject to} \quad & \mathbf{A}\mathbf{x} \preceq \mathbf{b}, \label{p1}
   \end{align}
with $\mathbf{A} \in \mathbb{R}^{B \times {N}M}$, $\mathbf{b} \in \mathbb{R}^{B}$.

\begin{Assumption}\label{a0}
The edges in ${G}$ and $G_i$ are undirected, and both ${G}$ and $G_i$ are connected, $\forall i\in V$.
\end{Assumption}

\begin{Assumption}\label{a1}
Both $f_{ij}: \mathbb{R}^M \rightarrow (-\infty,+\infty]$ and $g_{ij}: \mathbb{R}^M \rightarrow (-\infty,+\infty]$ are proper, convex, and closed extended real-valued functions; $f_{ij}$ is $\sigma_{ij}$-strongly convex, $\sigma_{ij}>0$, $\forall i\in {V}, j\in V_i$.
\end{Assumption}

\begin{Assumption}\label{a3a}
(Constraint Qualification) There exists an $\breve{\mathbf{x}} \in \mathbf{relint} S$ such that $\mathbf{A}\breve{\mathbf{x}} \preceq \mathbf{b}$, where $S \triangleq \prod_{i\in V} \bigcap_{j \in V_i} S_{ij} $ with $S_{ij}$ being the domain of $f_{ij} + g_{ij}$.
\end{Assumption}

\begin{Remark}\label{r0}
The composite cost function $f_{ij} + g_{ij}$ is a generalization of many cost functions in practical problems  \cite{beck2014fast,hans2009bayesian,zhao2017scope}. For instance, we can consider a local feasible region $X_{ij} \subseteq \mathbb{R}^M$ for the $j$th agent in cluster $i$, where $X_{ij}$ is non-empty, convex, and closed.
Then we can let $g_{ij}$ be an indicator function $\mathbb{I}_{X_{ij}}$, where $\mathbb{I}_{X_{ij}} (\mathbf{x}_{i}) \triangleq \left\{\begin{array}{ll}
                    0 & \hbox{if $\mathbf{x}_{i} \in X_{ij}$} \\
                    +\infty  & \hbox{otherwise}
                  \end{array}
                  \right.$\cite{notarnicola2016asynchronous}.
\end{Remark}

\begin{Remark}\label{r0-1}
For the comparison purpose, we consider a conventional composite DOP
\begin{align}
  \mathrm{(P1+)} \quad  \min\limits_{\mathbf{x}} \sum_{i \in {V}} (f_{i}(\mathbf{x}_{i}) + g_{i} (\mathbf{x}_{i})) \quad
   \hbox{subject to} \quad  \mathbf{A}\mathbf{x} \preceq \mathbf{b}, \nonumber
   \end{align}
where $f_i+g_i$ is the composite cost function of agent $i$. Problem (P1) can be viewed as a generalization of Problem (P1+) by extending the agent who manages $\mathbf{x}_i$ to cluster $i$ with multiple agents. In case there is only one agent in each cluster, Problem (P1) is equivalent to Problem (P1+). To optimize the cost function of Problem (P1), the agents in each cluster need to generate a consensual decision, {e.g.}, $\mathbf{x}_i$, which is involved in the coupling constraint $\mathbf{A}\mathbf{x} \preceq \mathbf{b}$. Therefore, the interactions among the agents within each cluster and across different clusters should be considered simultaneously.
\end{Remark}

\section{Optimization Algorithm Development}\label{sa4}

In this section, we propose an Asyn-DDPG algorithm for solving the problem of interest.

\subsection{Dual Problem}\label{sec1}

To realize distributed computations, we decouple the variable of clusters by defining $\mathbf{y}_{ij} \in \mathbb{R}^M$ as the estimate of $\mathbf{x}_i$ by the $j$th agent in cluster $i$. Then the collection of the estimates in cluster $i$ can be $\mathbf{y}_i \triangleq [\mathbf{y}^{\top}_{i1},...,\mathbf{y}^{\top}_{in_i}]^{\top} \in \mathbb{R}^{n_iM}$ and the collection of the overall estimates can be $\mathbf{y}  \triangleq [\mathbf{y}^{\top}_1,...,\mathbf{y}^{\top}_N]^{\top} \in \mathbb{R}^{\sum_{i\in V}n_iM}$.
To achieve the solution to Problem (P1), the estimates in each cluster should reach a consensus. Then by the consensus protocol in cluster $i$: $\mathbf{x}_i = \mathbf{y}_{i1} = ... = \mathbf{y}_{in_i}$, we have $\mathbf{A}\mathbf{x} = \sum_{i \in V} \mathbf{A}_i {\mathbf{x}}_i = \sum_{i \in V} \frac{\mathbf{1}^{\top}_{n_i} \otimes \mathbf{A}_i}{n_i} \mathbf{y}_i =  \sum_{i \in V} {\mathsf{A}_i\mathbf{y}_i} =  {\mathsf{A}\mathbf{y}}$,
where $\mathsf{A}_i \triangleq \frac{\mathbf{1}^{\top}_{n_i} \otimes \mathbf{A}_i}{n_i} \in \mathbb{R}^{B \times n_iM}$ and $\mathsf{A} \triangleq [{\mathsf{A}_1},...,{\mathsf{A}_N}] \in \mathbb{R}^{B \times \sum_{i \in V}n_i M}$ with $\mathbf{A}_i \in \mathbb{R}^{B \times M}$ being the $i$th column block of $\mathbf{A}$ ({i.e.}, $\mathbf{A}=[\mathbf{A}_1,...,\mathbf{A}_i,...,\mathbf{A}_N]$).

Note that the consensus constraint of $\mathbf{y}_{ij}$ in cluster $i$ can be equivalently written as $\mathsf{L}^i\mathbf{y}_i=\mathbf{0}$, where $\mathsf{L}^i \triangleq \mathbf{L}^i \otimes \mathbf{I}_M \in \mathbb{R}^{n_iM \times n_iM}$. Then the consensus-based optimization problem of the network $G$ can be formulated as
\begin{align}
  \mathrm{(P2)} \quad  \min\limits_{\mathbf{y},\mathbf{e}} \quad  & \sum_{i \in {V}} \sum_{j \in {V}_i} (f_{ij}(\mathbf{y}_{ij}) + g_{ij} (\mathbf{e}_{ij})) \nonumber \\
   \hbox{subject to} \quad & \mathbf{y}_{ij} = \mathbf{e}_{ij}, \quad \forall i \in V, j \in V_i,  \\
   & \mathsf{L}^i\mathbf{y}_i=\mathbf{0},\quad \forall i \in V,\\
   & \mathsf{A} \mathbf{y} \preceq \mathbf{b}, \label{ab1}
   \end{align}
where $\mathbf{e}_{ij} \in \mathbb{R}^{M}$ is a slack variable, $\mathbf{e}_{i} \triangleq [\mathbf{e}^{\top}_{i1},..., \mathbf{e}^{\top}_{in_i}]^{\top}$, and $\mathbf{e} \triangleq [\mathbf{e}^{\top}_{1},...,\mathbf{e}^{\top}_{N}]^{\top}$. Essentially, (\ref{ab1}) reconstructs (\ref{p1}) with variable $\mathbf{y}$ without affecting the nature of the constraint when certain consensus is achieved in each cluster.

The Lagrangian function of Problem (P2) can be given by $ \mathcal{L}_F(\mathbf{y}, \mathbf{e}, \bm{\mu}, \bm{\nu},\bm{\phi})
\triangleq \sum_{i \in {V}} (\sum_{j \in {V}_i} $ $(f_{ij}(\mathbf{y}_{ij}) + g_{ij}(\mathbf{e}_{ij}) +  \bm{\mu}_{ij}^{\top} (\mathbf{y}_{ij} - \mathbf{e}_{ij})) + \bm{\nu}^{\top}_{i}\mathsf{L}^i \mathbf{y}_i )+ \bm{\phi}^{\top} (\mathsf{A}\mathbf{y} - \mathbf{b})
= \sum_{i \in {V}} (\sum_{j \in {V}_i} (f_{ij}(\mathbf{y}_{ij}) + g_{ij}(\mathbf{e}_{ij}) + \bm{\mu}_{ij}^{\top} (\mathbf{y}_{ij} - \mathbf{e}_{ij})) + \bm{\nu}^{\top}_{i}\mathsf{L}^i \mathbf{y}_i + \bm{\phi}^{\top} \mathsf{A}_i\mathbf{y}_i ) - \bm{\phi}^{\top}  \mathbf{b} =  \sum_{i \in {V}} \sum_{j \in {V}_i} (f_{ij}(\mathbf{y}_{ij}) + g_{ij}(\mathbf{e}_{ij}) + \bm{\mu}_{ij}^{\top} (\mathbf{y}_{ij} - \mathbf{e}_{ij} ) + \bm{\nu}^{\top}_{i}\mathsf{L}^i_{j} \mathbf{y}_{ij} + \bm{\phi}^{\top} \mathsf{A}_{ij}\mathbf{y}_{ij}  - \bm{\phi}^{\top}\mathbf{b}_{ij})$, where $\bm{\mu}_{ij} \in \mathbb{R}^M$ and $\bm{\nu}_i \in \mathbb{R}^{n_iM}$ and $\bm{\phi} \in \mathbb{R}^{B}_+$ are Lagrangian multipliers,
$\bm{\mu}_{i} \triangleq [\bm{\mu}^{\top}_{i1},..., \bm{\mu}^{\top}_{in_i}]^{\top}, \bm{\mu} \triangleq [\bm{\mu}^{\top}_{1},...,\bm{\mu}^{\top}_{N}]^{\top}, \bm{\nu} \triangleq [\bm{\nu}^{\top}_{1},...,\bm{\nu}^{\top}_{N}]^{\top}$,
$\mathsf{A}_{ij} \in \mathbb{R}^{B \times M}$ is the $j$th column block of $\mathsf{A}_{i}$ ({i.e.}, $\mathsf{A}_{i} = [\mathsf{A}_{i1},...,\mathsf{A}_{ij},...,\mathsf{A}_{in_i}]$), $\mathsf{L}^i_{j} \in \mathbb{R}^{n_iM \times M}$ is the $j$th column block of $\mathsf{L}^i$ ({i.e.}, $\mathsf{L}^i = [\mathsf{L}^i_{1},...,\mathsf{L}^i_{j},...,\mathsf{L}^i_{n_i}]$), $\sum_{i\in V} \sum_{j\in V_i}\mathbf{b}_{ij} = \mathbf{b}$.

Then the dual function can be obtained by $D(\bm{\mu},\bm{\nu},\bm{\phi})
\triangleq  \min\limits_{\mathbf{y}, \mathbf{e}} \sum_{i \in {V}} \sum_{j \in {V}_i} (f_{ij}(\mathbf{y}_{ij}) + g_{ij}(\mathbf{e}_{ij})  + \bm{\mu}_{ij}^{\top} (\mathbf{y}_{ij} - \mathbf{e}_{ij} ) + \bm{\nu}^{\top}_{i}\mathsf{L}^i_{j} \mathbf{y}_{ij} + \bm{\phi}^{\top} \mathsf{A}_{ij}\mathbf{y}_{ij}  - \bm{\phi}^{\top} \mathbf{b}_{ij} )
=  \min\limits_{\mathbf{y}, \mathbf{e}} \sum_{i \in {V}} \sum_{j \in {V}_i} (f_{ij}(\mathbf{y}_{ij}) + (\bm{\mu}_{ij} + \mathsf{L}^{i\top}_{j}\bm{\nu}_{i} + \mathsf{A}^{\top}_{ij} \bm{\phi})^{\top} \mathbf{y}_{ij} - \mathbf{b}_{ij}^{\top} \bm{\phi} + g_{ij}(\mathbf{e}_{ij}) - \bm{\mu}^{\top}_{ij}\mathbf{e}_{ij})
= \sum_{i \in {V}} \sum_{j \in {V}_i} (- f^{\circ}_{ij}(-\bm{\mu}_{ij} - \mathsf{L}^{i\top}_{j}\bm{\nu}_{i} - \mathsf{A}^{\top}_{ij} \bm{\phi}) - \mathbf{b}_{ij}^{\top} \bm{\phi} - g_{ij}^{\circ}(\bm{\mu}_{ij}))$,
where $f^{\circ}_{ij}$ and $g^{\circ}_{ij}$ are the Fenchel conjugates of $f_{ij}$ and $g_{ij}$, respectively. Hence, the dual problem of Problem (P2) is
\begin{align}
\mathrm{(P3)} \quad \max\limits_{\bm{\mu}, \bm{\nu}, \bm{\phi} \succeq \mathbf{0}} & D(\bm{\mu},\bm{\nu},\bm{\phi}). \nonumber
\end{align}
{Let $(\bm{\mu}^*,\bm{\nu}^*,\bm{\phi}^*)$ be an optimal solution to Problem (P3). In addition, let ${Y_i} \subset \mathbb{R}^{n_iM}$ and $J \subset \mathbb{R}^B_+$ be two non-empty, convex, and compact zones that estimated by practitioners, such that $\bm{\nu}_i^* \in {Y_i}$ and $\bm{\phi}^* \in J$, $\forall i \in V$. The range of ${Y_i}$ and $J$ can be larger or smaller depending on how the practitioners are confident about the location of $\bm{\nu}_i^*$ and $\bm{\phi}^*$. Then there exist some real values $\iota_Y \triangleq  \max_{i \in {V}} \{ \max_{\bm{\nu}_i \in {Y_i}}\| \bm{\nu}_{i} \| \}$ and $\iota_J \triangleq  \max_{\bm{\phi} \in J } \| \bm{\phi}\|$. Note that the estimated zones exist as long as $\bm{\nu}_i^*$ and $\bm{\phi}^*$ exist without any additional assumption on the parameters of the primal problem. In case $\bm{\nu}_i^* \notin {Y_i}$ and/or $\bm{\phi}^* \notin J$ (i.e., the estimation of practitioners is not accurate enough), there may be some gaps in the final optimization result, which is out of the scope of this research work.}

\begin{Remark}
{The estimation of ${Y_i}$ and $J$ serves as an intermediate step when formulating the problem of interest to facilitate the algorithm development.} Similar settlements can be referred to in \cite{notarnicola2019constraint,notarnicola2017duality}. In practice, the estimation of ${Y_i} $ and ${J}$ relies on the experience in specific problems. For example, in some social welfare optimization problems in the electricity market, the optimal dual variables can be the settled energy prices \cite{samadi2010optimal}, whose range can be estimated easily with historical prices as long as the parameters of the market do not vary seriously over the years.
\end{Remark}

Based on the estimated zones, Problem (P3) is equivalent to

\begin{align}
\mathrm{(P4)} \quad \min\limits_{\bm{\mu}, \bm{\nu}, \bm{\phi} } & \sum_{i \in {V}} \sum_{j \in {V}_i} \big(f^{\circ}_{ij}(-\bm{\mu}_{ij} - \mathsf{L}^{i\top}_{j}\bm{\nu}_{i} - \mathsf{A}^{\top}_{ij} \bm{\phi}) + \mathbf{b}_{ij}^{\top} \bm{\phi} \big)  \nonumber \\
& + \sum_{i \in {V}} \sum_{j \in {V}_i} (g_{ij}^{\circ}(\bm{\mu}_{ij}) + \mathbb{I}_{Y_i}(\bm{\nu}_{i})+ \mathbb{I}_{J}(\bm{\phi})). \nonumber
\end{align}
In Problem (P4), we introduce indicator functions $\mathbb{I}_{Y_i}$ and $\mathbb{I}_{J}$ without affecting the solution compared with (P3) with $\bm{\nu}_i^* \in {Y_i}$ and $\bm{\phi}^* \in J$, $\forall i \in V$. The constraint $\bm{\phi} \succeq \mathbf{0}$ in (P3) is addressed by $\mathbb{I}_{J}$ since $J \subset \mathbb{R}^B_+$.

Then to solve Problem (P4) in a distributed manner, we introduce the local estimates for handling the common variables. Specifically, let $\bm{\alpha}_{ij} \triangleq [\bm{\mu}_{ij}^{\top},\bm{\gamma}^{\top}_{ij},\bm{\theta}_{ij}^{\top}]^{\top}$ be maintained by the $j$th agent in cluster $i$, where $\bm{\gamma}_{ij}$ and $\bm{\theta}_{ij}$ are the local estimates of $\bm{\nu}_i$ and $\bm{\phi}$, respectively.
Then with the given ${Y_i}$ and $J$, we have $\iota_Y =  \max_{i \in {V}} \{ \max_{\bm{\gamma}_{ij} \in {Y_i}}\| \bm{\gamma}_{ij} \| \}$ and $\iota_J =  \max_{\bm{\theta}_{ij} \in J } \| \bm{\theta}_{ij} \|$.
In addition, define $ \bm{\gamma}_{i} \triangleq [\bm{\gamma}^{\top}_{i1},\bm{\gamma}^{\top}_{i2},..., \bm{\gamma}^{\top}_{in_i}]^{\top}$, $ \bm{\gamma} \triangleq [\bm{\gamma}^{\top}_{1},..., \bm{\gamma}^{\top}_{N}]^{\top}, \bm{\theta}_{i} \triangleq [\bm{\theta}^{\top}_{i1},\bm{\theta}^{\top}_{i2},..., \bm{\theta}^{\top}_{in_i}]^{\top}$, $ \bm{\theta} \triangleq [\bm{\theta}^{\top}_{1},..., \bm{\theta}^{\top}_{N}]^{\top}, \bm{\alpha}_i \triangleq [\bm{\alpha}^{\top}_{i1} ,\bm{\alpha}^{\top}_{i2} ,...,\bm{\alpha}^{\top}_{in_i}]^{\top}$, $ \bm{\alpha} \triangleq [\bm{\alpha}^{\top}_1,...,\bm{\alpha}^{\top}_{N}]^{\top}, \mathbf{N}_{ij} \triangleq [\mathbf{O}_{B \times (M+n_iM)}, \mathbf{I}_B]$, $ \mathbf{N}_i \triangleq \mathbf{I}_{n_i} \otimes \mathbf{N}_{ij}, \mathbf{M}_{ij} \triangleq [\mathbf{O}_{n_iM \times M},\mathbf{I}_{n_iM},\mathbf{O}_{n_iM \times B}]$, $ \mathbf{M}_i \triangleq \mathbf{I}_{n_i} \otimes \mathbf{M}_{ij},
\mathbf{N} \triangleq \mathrm{diag}[\mathbf{N}_{1},...,\mathbf{N}_{N}]$, $\mathbf{M} \triangleq \mathrm{diag}[\mathbf{M}_{1},...,\mathbf{M}_{N}],
\mathbf{W}_{ij} \triangleq [-\mathbf{I}_M,-\mathsf{L}^{i\top}_{j}, -\mathsf{A}^{\top}_{ij}]$, $ \mathbf{D}_{ij} \triangleq [\mathbf{0}^{\top}_M,\mathbf{0}^{\top}_{n_iM}, \mathbf{b}_{ij}^{\top}]$, $ \mathbf{R}_{ij} \triangleq [\mathbf{I}_M,\mathbf{O}_{M \times n_iM}, \mathbf{O}_{M \times B}]$.

With the local estimates, we have {cluster-based} cost functions and consensus protocol. In particular, based on Problem (P4), the cost function with the local estimates of agent $n_{ij}$ (relabeled index) can be given by $s_{ij}(\bm{\alpha}_{ij}) + r_{ij}(\bm{\alpha}_{ij})$, where $s_{ij}(\bm{\alpha}_{ij})  \triangleq f^{\circ}_{ij}(\mathbf{W}_{ij} \bm{\alpha}_{ij}) +\mathbf{D}_{ij} \bm{\alpha}_{ij}$ and $ r_{ij}(\bm{\alpha}_{ij})  \triangleq g_{ij}^{\circ}(\mathbf{R}_{ij} \bm{\alpha}_{ij}) + \mathbb{I}_{Y_i}(\mathbf{M}_{ij}\bm{\alpha}_{ij}) + \mathbb{I}_{J}(\mathbf{N}_{ij}\bm{\alpha}_{ij})$.
Note that $s_{ij}$ and $r_{ij}$ are convex by the definition of Fenchel conjugate \cite{beck2014fast}. Then the global cost function can be defined by $H ( {\bm{\alpha}}) \triangleq H_s( {\bm{\alpha}}) + H_r ( {\bm{\alpha}})$, where $H_s( {\bm{\alpha}})  \triangleq  \sum_{i \in {V}} \sum_{j \in {V}_i} s_{ij}(\bm{\alpha}_{ij})$ and $H_r( {\bm{\alpha}})  \triangleq \sum_{i \in {V}} \sum_{j \in {V}_i} r_{ij}(\bm{\alpha}_{ij})$.
In addition, the consensus of ${\bm{\gamma}}_{ij}$ in $V_i$ and ${\bm{\theta}}_{ij}$ in $\widehat{V}$ can be characterized by
\begin{align}\label{}
& \bm{\gamma}_{ij} - \bm{\gamma}_{il} =\mathbf{0}, \quad \forall l \in \Omega_{ij}, i \in V, j \in V_i, \label{a2} \\
& {\bm{\theta}}_{ij} - {\bm{\theta}}_{uv} = \mathbf{0}, \quad  \forall n_{uv} \in \widehat{\Omega}_{ij}, i \in V, j \in V_i, \label{a3}
\end{align}
respectively, where $\Omega_{ij} \triangleq \{l| (j,l) \in E_i, l >j\}$ and $\widehat{\Omega}_{ij} \triangleq \{k| (n_{ij},k) \in E, k> n_{ij} \}$. (\ref{a2}) means that the consensus between the $l$th agent in cluster $i$ and the $j$th agent in the same cluster is achieved. Similarly, the consensus between the $n_{uv}$th agent (relabeled index) and the $n_{ij}$th agent (relabeled index) in $\widehat{V}$ is achieved by (\ref{a3}). $\Omega_{ij}$ and $\widehat{\Omega}_{ij}$ are empty if no consensus constraints exist in (\ref{a2}) and (\ref{a3}), respectively. {Based on the discussions above, the cluster-based optimization problem (P4) with the local estimate of the agents can be reformulated as
\begin{align}\label{}
\mathrm{(P5)} \quad \min_{\bm{\alpha}} H(\bm{\alpha}) \quad \hbox{subject to} \quad (\ref{a2})\hbox{ and }(\ref{a3}). \nonumber
\end{align}}



Note that $\bm{\gamma}_{ij} = \mathbf{M}_{ij}\bm{\alpha}_{ij}$, $\bm{\theta}_{ij} =  \mathbf{N}_{ij}\bm{\alpha}_{ij}$, $\bm{\gamma}_i = \mathbf{M}_i\bm{\alpha}_i$, $\bm{\gamma} = \mathbf{M}\bm{\alpha}$, $\bm{\theta}_i =  \mathbf{N}_i\bm{\alpha}_i$, and $\bm{\theta} =  \mathbf{N}\bm{\alpha}$. Then (\ref{a2}) and (\ref{a3}) can be written in compact forms with the help of the incidence matrix \cite{dimarogonas2010stability}. By defining $\mathsf{G}_i \triangleq \mathbf{G}_i^{\top} \otimes \mathbf{I}_{n_iM}$, (\ref{a2}) can be represented by $\mathsf{G}_i \bm{\gamma}_{i}= \mathsf{G}_i \mathbf{M}_i \bm{\alpha}_{i} = \mathbf{0}$, $\forall i \in V$. In addition, one can construct $\mathsf{G}\bm{\gamma} = \mathsf{G} \mathbf{M} \bm{\alpha} = \mathbf{0}$ by including all clusters where $\mathsf{G} \triangleq \mathrm{diag}[\mathsf{G}_1,...,\mathsf{G}_{N}]$. Similarly, (\ref{a3}) can be represented by $ \widehat{\mathsf{G}} {\bm{\theta}}= \widehat{\mathsf{G}}\mathbf{N}{\bm{\alpha}} = \mathbf{0}$, where $\widehat{\mathsf{G}} \triangleq \widehat{\mathbf{G}}^{\top} \otimes \mathbf{I}_B$. In addition, let $\mathsf{G}_i^j$ be the row sub-block of $\mathsf{G}$, which contains the edge information in cluster $i$, associated with which $j$ is the smaller vertex index. Similarly, $\widehat{\mathsf{G}}^j_{i}$ is the row sub-block of $\widehat{\mathsf{G}}$, which contains the edge information associated with which $n_{ij}$ is the smaller relabeled vertex index. The row dimension of $\mathsf{G}_i^j$ and $\widehat{\mathsf{G}}_i^j$ is zero if $\Omega_{ij} = \emptyset$ and $\widehat{\Omega}_{ij} = \emptyset$, respectively.
To facilitate the following analysis, we define $\mathbf{Z} \triangleq \left[
  \begin{array}{l}
  \mathsf{G}\mathbf{M} \\
  \widehat{\mathsf{G}} \mathbf{N} \\
  \end{array}
\right]$.
Then (\ref{a2}) and (\ref{a3}) can be jointly represented by $\mathbf{Z}\bm{\alpha} = \left[
                          \begin{array}{c}
                            \mathsf{G}\bm{\gamma} \\
                            \widehat{\mathsf{G}} {\bm{\theta}} \\
                          \end{array}
                        \right]
=\mathbf{0}$.
Define $\mathbf{Z}_{ij} \triangleq
\left[
  \begin{array}{c}
    \mathsf{G}_i^{j} \mathbf{M}\\
    \widehat{\mathsf{G}}_i^{j} \mathbf{N}\\
  \end{array}
\right]$, $\forall i \in V$, $j \in V_i$. Then we have $\mathbf{Z}_{ij} \bm{\alpha} = \left[
  \begin{array}{c}
    \mathsf{G}_i^{j} \mathbf{M} \bm{\alpha}\\
    \widehat{\mathsf{G}}_i^{j} \mathbf{N} \bm{\alpha} \\
  \end{array}
\right] =
\left[
  \begin{array}{c}
    \mathsf{G}_i^{j} \bm{\gamma} \\
    \widehat{\mathsf{G}}_i^{j} \bm{\theta} \\
  \end{array}
\right]$.
Due to the boundedness of $\bm{\gamma}$ and ${\bm{\theta}}$, the upper bound of $ \parallel \mathbf{Z}_{ij} \bm{\alpha} \parallel$ can be obtained by $\max \{\parallel \mathbf{Z}_{ij} \bm{\alpha} \parallel\}_{i\in V, j\in V_i}
=  \max \bigg\{ \sqrt{\| \mathsf{G}_i^{j} \bm{\gamma} \|^2 + \| \widehat{\mathsf{G}}_i^{j} {\bm{\theta}} \|^2} \bigg\}_{i\in V, j\in V_i}
\leq  \max \bigg\{ \sqrt{|\widehat{V} | \| \mathsf{G}_i^{j} \|_2^2 \iota^2_Y + |\widehat{V} | \| \widehat{\mathsf{G}}_i^{j} \|_2^2 \iota^2_J} \bigg\}_{i\in V, j\in V_i}
\triangleq \iota$.

To facilitate the algorithm development, we introduce the following auxiliary variables. $ \bm{\xi}_{ijl} \in \mathbb{R}^{n_iM}$, $ \bm{\xi}_{ij} \triangleq [ \bm{\xi}^{\top}_{ijl_1},...,\bm{\xi}^{\top}_{ijl_{|\Omega_{ij}|}}]^{\top}$, $ l_{(\cdot)} \in \Omega_{ij}$, $\bm{\zeta}_{ijk} \in \mathbb{R}^B$, $\bm{\zeta}_{ij} \triangleq [\bm{\zeta}^{\top}_{ijk_{1}},...,\bm{\zeta}^{\top}_{ijk_{|\widehat{\Omega}_{ij}|}}]^{\top}, k_{(\cdot)} \in \widehat{\Omega}_{ij}$, $\bm{\xi}_{i} \triangleq [\bm{\xi}^{\top}_{i1},\bm{\xi}^{\top}_{i2},...,\bm{\xi}^{\top}_{in_i}]^{\top}$, $ \bm{\xi} \triangleq [\bm{\xi}^{\top}_{1},...,\bm{\xi}^{\top}_{N}]^{\top}, \bm{\zeta}_{i} \triangleq [\bm{\zeta}^{\top}_{i1},\bm{\zeta}^{\top}_{i2},...,\bm{\zeta}^{\top}_{in_i}]^{\top}$, $\bm{\zeta} \triangleq [\bm{\zeta}^{\top}_{1},...,\bm{\zeta}^{\top}_{N}]^{\top}, \bm{\omega} \triangleq [\bm{\xi}^{\top},\bm{\zeta}^{\top}]^{\top}$.
Here, the sequence $\{l_1,...,l_{|\Omega_{ij}|}\}$ is decided as follows. $\forall l_m,l_n \in \Omega_{ij}$, if $l_m>l_n$, then $m>n$, and vice versa. Similarly, in $\{k_{1},...,k_{|\widehat{\Omega}_{ij}|}\}$, $\forall k_{m},k_{n} \in \widehat{\Omega}_{ij}$, if $k_{m}>k_{n}$, then $m>n$, and vice versa. $\bm{\xi}_{ij}$ and $\bm{\zeta}_{ij}$ do not exist if $\Omega_{ij}$ and $\widehat{\Omega}_{ij}$ are empty, respectively.



\subsection{Asyn-DDPG Algorithm Development}

In the following, we develop an Asyn-DDPG algorithm for solving Problem (P5) by considering the communication delays between the neighboring agents.
%
In this work, the delay of the state of agent $n_{uv}$ (relabeled index) accessed by agent $n_{ij} \in \widehat{V}^{n_{uv}}$ at step $t$ is defined by $q_{uv}^{ij,t} \in \mathbb{N}$, which satisfies the following assumption.

\begin{Assumption}\label{as5}
The communication delay between any two neighboring agents is upper bounded by $q \triangleq \max\{ q_{uv}^{ij,t}\}_{n_{ij} \in \widehat{V},n_{uv} \in \widehat{V}^{n_{ij}}, t \in \mathbb{N}}$. The delayed states will reach each agent in a sequential order in which they were sent.
\end{Assumption}

\begin{Remark}
The assumption on the bounded delays is widely discussed in the existing works \cite{chang2016asynchronous,wang2022composite,tian2020achieving,zhou2018distributed}. 
To access the outdated data, buffers can be established \cite{zhang2019asyspa,feyzmahdavian2016asynchronous,tian2020achieving}, where the data is stamped with step index $t$, such that the historical data at a specific step can be read by the agents.
\end{Remark}

Then we define some {outdated} versions of variables. $ \underline{\bm{\alpha}}^t_{ij} = [
    \underline{\bm{\mu}}_{ij}^{t\top},
    \underline{\bm{\gamma}}_{ij}^{t\top},
    \underline{\bm{\theta}}_{ij}^{t\top} ]^{\top} \triangleq {\bm{\alpha}}_{ij}^{(t-d)^+} =[ \bm{\mu}_{ij}^{(t-d)^+\top},
    \bm{\gamma}_{ij}^{(t-d)^+\top},
    {\bm{\theta}}_{ij}^{(t-d)^+\top} ]^{\top}$, $ \underline{\bm{\xi}}^t_{ijl} \triangleq {\bm{\xi}}_{ijl}^{(t-d)^+}, \underline{\bm{\zeta}}^t_{ijk} \triangleq {\bm{\zeta}}_{ijk}^{(t-d)^+}$, where $d \triangleq 2q+ 1$, $i \in V, j\in V_i, l\in \Omega_{ij}, k \in \widehat{\Omega}_{ij}$. Here, we use notation ``$\underline{\  }$'' to represent an outdated state generated $d$ steps earlier than step $t$ but no earlier than step $0$.
Then we define the following vectors. $ \underline{\bm{\alpha}}_i^t \triangleq [\underline{\bm{\alpha}}^{t\top}_{i1},\underline{{\bm{\alpha}}}^{t\top}_{i2},...,\underline{{\bm{\alpha}}}^{t\top}_{in_i}]^{\top} $, $\underline{\bm{\alpha}}^t \triangleq [\underline{\bm{\alpha}}^{t\top}_{1},..., \underline{\bm{\alpha}}^{t\top}_{N}]^{\top}$, $\underline{\bm{\xi}}^t_{ij} \triangleq [ \underline{\bm{\xi}}^{t\top}_{ijl_1},...,\underline{\bm{\xi}}^{t\top}_{ijl_{|\Omega_{ij}|}} ]^{\top}$, $\underline{\bm{\zeta}}^t_{ij} \triangleq [\underline{\bm{\zeta}}^{t\top}_{ijk_{1}},...,\underline{\bm{\zeta}}^{t\top}_{ijk_{|\widehat{\Omega}_{ij}|}} ]^{\top}$, $\underline{\bm{\xi}}^t_{i} \triangleq [\underline{\bm{\xi}}^{t\top}_{i1},...,\underline{\bm{\xi}}^{t\top}_{in_i} ]^{\top}$, $\underline{\bm{\xi}}^t \triangleq [\underline{\bm{\xi}}^{t\top}_{1},...,\underline{\bm{\xi}}^{t\top}_{N}]^{\top}$, $ \underline{\bm{\zeta}}^t_{i} \triangleq [\underline{\bm{\zeta}}^{t\top}_{i1},...,\underline{\bm{\zeta}}^{t\top}_{in_i}]^{\top}$, $ \underline{\bm{\zeta}}^t \triangleq [\underline{\bm{\zeta}}^{t\top}_{1},...,\underline{\bm{\zeta}}^{t\top}_{N}]^{\top}$, $ \underline{\bm{\omega}}^t \triangleq [\underline{\bm{\xi}}^{t\top},\underline{\bm{\zeta}}^{t\top}]^{\top}$.
\begin{algorithm}
\caption{\textbf{Asyn-DDPG algorithm}}\label{ax1}
\begin{algorithmic}[1]
\State \textbf{Initialize} $\bm{\alpha}^{0}$, $\bm{\omega}^{0}$, and step sizes $c_{ij}>0$ and $\pi_{ij}>0$, $\forall i \in {V}$, $j \in V_i$.
\For {$t= 0,1,2,...$}
\For {$i= 1,2,...,{N}$} // {\texttt{In parallel}}
\For {$j=1,2,...,n_i$} // {\texttt{{In parallel}}}
\State \textbf{\texttt{{UPDATE $\bm{\alpha}$}}}:
\State \textbf{Calculate} // \texttt{{With local state}}
\begin{align}\label{}
\mathbf{m}_{\bm{\mu}_{ij}}^t = & - \arg \max\limits_{\mathbf{n}} ( (\mathbf{W}_{ij} \bm{\alpha}^t_{ij})^{\top} \mathbf{n} - f_{ij}(\mathbf{n})), \label{m1} \\
\mathbf{m}_{\bm{\gamma}_{ij}}^t =  & -\mathsf{L}^{i}_{j} \arg \max\limits_{\mathbf{n}} ( (\mathbf{W}_{ij} \bm{\alpha}^t_{ij})^{\top} \mathbf{n} - f_{ij}(\mathbf{n})),  \label{m2} \\
\mathbf{m}_{\bm{\theta}_{ij}}^t = & -\mathsf{A}_{ij} \arg \max\limits_{\mathbf{n}} ( (\mathbf{W}_{ij} \bm{\alpha}^t_{ij})^{\top} \mathbf{n} - f_{ij}(\mathbf{n})) \nonumber \\
& +  \mathbf{b}_{ij} .  \label{m3}
\end{align}
\State Upon the arrival of $\underline{\bm{\xi}}_{ijl}^t$, $\underline{\bm{\xi}}_{il^{'}j}^t$, $\underline{\bm{\gamma}}_{il}^t$, $\underline{\bm{\gamma}}_{il^{'}}^t$,  ${\underline{\bm{\zeta}}}_{ijn_{uv}}^t$, ${\underline{\bm{\zeta}}}_{u'v'n_{ij}}^t$, $\underline{{\bm{\theta}}}_{uv}^t$, and $\underline{{\bm{\theta}}}_{u'v'}^t$, $\forall l \in \Omega_{ij}$, $l^{'} \in \Omega^{\sharp}_{ij}$, $n_{uv} \in \widehat{\Omega}_{ij},n_{u'v'} \in \widehat{\Omega}^{\sharp}_{ij}$, \textbf{update} // \texttt{{With outdated state of neighbors}}
\begin{align}
 \bm{\mu}_{ij}^{t+1}  = & \mathbf{prox}^{c_{ij}}_{g_{ij}^{\circ}} [\bm{\mu}_{ij}^t - c_{ij} \mathbf{m}_{\bm{\mu}_{ij}}^t ], \label{1} \\
 \bm{\gamma}_{ij}^{t+1}  = & \mathcal{P}_{Y_i} \bigg[\bm{\gamma}_{ij}^t - c_{ij}  \bigg( \mathbf{m}_{\bm{\gamma}_{ij}}^t +  \sum_{ l \in \Omega_{ij}} \underline{\bm{\xi}}_{ijl}^t  - \sum_{ l^{'} \in \Omega^{\sharp}_{ij}} \underline{\bm{\xi}}_{il^{'}j}^t \nonumber \\
 &  + \pi_{ij} \sum_{l \in \Omega_{ij}}  ( \underline{\bm{\gamma}}_{ij}^t - \underline{\bm{\gamma}}_{il}^t)  + \sum_{l^{'} \in \Omega^{\sharp}_{ij}} \pi_{il^{'}}  (\underline{\bm{\gamma}}_{ij}^t - \underline{\bm{\gamma}}_{il^{'}}^t ) \bigg) \bigg] , \label{2}  \\
 \bm{\theta}_{ij}^{t+1} = & \mathcal{P}_{J} \bigg[\bm{\theta}_{ij}^t - c_{ij}  \bigg( \mathbf{m}_{\bm{\theta}_{ij}}^t  + \sum_{n_{u v} \in \widehat{\Omega}_{ij}}
{\underline{\bm{\zeta}}}_{ijn_{u v}}^t \nonumber \\
 & \quad \quad  -  \sum_{n_{u'v'} \in \widehat{\Omega}^{\sharp}_{ij}} {\underline{\bm{\zeta}}}_{u'v'n_{ij}}^t  + \pi_{ij}\sum_{n_{uv} \in \widehat{\Omega}_{ij} } (\underline{{\bm{\theta}}}_{ij}^t - \underline{{\bm{\theta}}}_{u v}^t  )  \nonumber \\
 &  \quad \quad + \sum_{n_{u'v'} \in \widehat{\Omega}^{\sharp}_{ij}} {\pi}_{u'v'} (\underline{{\bm{\theta}}}_{ij}^t- \underline{{\bm{\theta}}}_{u'v'}^t) \bigg) \bigg]. \label{3}
\end{align}
\textbf{Broadcast} the result to the neighbors.
\State \textbf{\texttt{UPDATE $\bm{\omega}$}}:
\State Upon the arrival of $\bm{\gamma}^{t+1}_{il}$ and ${\bm{\theta}}^{t+1}_{uv}$, $\forall l \in \Omega_{ij}$, $n_{uv} \in \widehat{\Omega}_{ij}$, \textbf{update} // \texttt{{With outdated state of neighbors}}
\begin{align}\label{}
& \bm{\xi}_{ijl}^{t+1} =   \underline{\bm{\xi}}_{ijl}^t  + \pi_{ij} (\bm{\gamma}_{ij}^{t+1} - \bm{\gamma}_{il}^{t+1}), \label{4}   \\
&   \bm{\zeta}_{ijn_{uv}}^{t+1} =   \underline{\bm{\zeta}}_{ijn_{uv}}^t + \pi_{ij} ( {\bm{\theta}}_{ij}^{t+1} - {\bm{\theta}}_{uv}^{t+1}). \label{5}
\end{align}
\textbf{Broadcast} the result to the neighbors.
\EndFor
\EndFor
\EndFor
\State \textbf{Obtain} the outputs, denoted by $\bar{{\bm{\mu}}}_{ij},$ $\bar{{\bm{\gamma}}}_{ij}$, $\bar{{\bm{\theta}}}_{ij}$, $\bar{{\bm{\xi}}}_{ijl}$, and $\bar{{\bm{\zeta}}}_{ijn_{uv}}$, under certain convergence criterion, $\forall i \in V$, $j\in V_i$, $l\in \Omega_{ij}$, $n_{uv} \in \widehat{\Omega}_{ij}$.
\State \textbf{Calculate} the estimates of the solution to Problem (P1): $\bar{\mathbf{y}}_{ij}  = \arg \min_{\mathbf{n}} f_{ij}(\mathbf{n}) + (\bar{\bm{\mu}}_{ij} + \mathsf{L}^{i\top}_{j} \bar{\bm{\gamma}}_{ij} + \mathsf{A}^{\top}_{ij} \bar{\bm{\theta}}_{ij})^{\top} \mathbf{n}$, $\forall i \in V$, $j\in V_i$.
\end{algorithmic}
\end{algorithm}
To proceed, we define $\Omega^{\sharp}_{ij} \triangleq \{l| (j,l) \in E_i, l < j\}$ and $\widehat{\Omega}^{\sharp}_{ij} \triangleq \{k| (n_{ij},k) \in E, k < n_{ij} \}$, $\forall i \in {V}, j \in V_i$.
Then the proposed Asyn-DDPG algorithm for optimizing the problem of interest is designed in Algorithm \ref{ax1}. An illustrative data acquisition process of agent $n_{ij}$ (relabeled index) is shown in Fig. \ref{ccp}.



\begin{Remark}\label{a4}
The computation of (\ref{1}) can be further simplified by {extended Moreau decomposition} \cite[Thm. 6.45]{beck2017first}, which gives $ \bm{\mu}_{ij}^{t+1} =  \mathbf{prox}^{c_{ij}}_{g_{ij}^{\circ}} [\mathbf{s}^t_{ij} ]  =  \mathbf{s}^t_{ij}   - {c_{ij}}\mathbf{prox}^{{c^{-1}_{ij}}}_{g^{\circ \circ }_{ij}} [{c^{-1}_{ij}}\mathbf{s}^t_{ij} ]$, where $\mathbf{s}^t_{ij} \triangleq  \bm{\mu}_{ij}^t - c_{ij} \mathbf{m}_{\bm{\mu}_{ij}}^t$, $g^{\circ \circ }_{ij} = g_{ij} $ and $g^{\circ \circ }_{ij}$ is the biconjugate of $g_{ij}$ \cite[Sec. 3.3.2]{boyd2004convex}. In this case, the calculation of $g^{\circ}_{ij}$ is not required.
To compute $\bm{\alpha}^{t+1}$ and $\bm{\omega}^{t+1}$, we let $d=2q +1$, which is the upper bound of the sum of the following 3 values: 1) the number of instants for transmitting $\underline{\bm{\gamma}}^{t}$, $\underline{\bm{\theta}}^{t}$, $\underline{\bm{\xi}}^{t}$, and $\underline{\bm{\zeta}}^{t}$ among neighboring agents to compute $\bm{\alpha}^{t+1}$ based on (\ref{2}) and (\ref{3}) (upper bounded by $q$); 2) the number of instants for transmitting ${\bm{\gamma}}^{t+1}$ and ${\bm{\theta}}^{t+1}$ among neighboring agents to compute ${\bm{\omega}}^{t+1}$ based on (\ref{4}) and (\ref{5}) (upper bounded by $q$); 3) the number of instants of an update (equals 1).       Consequently, the computations on $\bm{\alpha}$ and $\bm{\omega}$ can be conducted successively along the time line since the information at the right-hand side of (\ref{1})-(\ref{5}) is available in the buffers for every $t$. {{In practice, the uniform delay bound can be negotiated among all the agents before applying the proposed algorithm since this value is time invariant.}}
\begin{Remark}
Existing consensus-based distributed optimization algorithms {without coupling constraints} are often based on a standard optimization scheme $\min\limits_{\mathbf{x}} \sum_{i \in V} f_i(\mathbf{x})$, where certain agreement on a common decision $\mathbf{x}$ is achieved  \cite{cao2021decentralized,sirb2018decentralized,tian2020achieving,tsianos2012distributed,zhang2019asyspa,chang2016asynchronous}. Differently in Problem (P4), the cost function is decided by agent $n_{ij}$'s (relabeled index) {private} variable $\bm{\mu}_{ij}$, {cluster-wide} common variable $\bm{\nu}_i$, and {network-wide} common variable $\bm{\phi}$, which is due to the {heterogeneous} objectives within and across the clusters. Then the consensus protocol is established among the agents in each cluster (w.r.t. $\bm{\nu}_i$) and the overall agents in the network (w.r.t. $\bm{\phi}$). This leads to a {partial} (with the existence of private $\bm{\mu}_{ij}$) and {mixed} (with both cluster-wide and network-wide agreements) consensus protocol, which can be more general than the conventional consensus schemes {only with} global agreements \cite{cao2021decentralized,sirb2018decentralized,tian2020achieving,tsianos2012distributed,zhang2019asyspa,chang2016asynchronous}. 
\end{Remark}
      %

\begin{figure}[H]
  \centering
  \includegraphics[width=10cm]{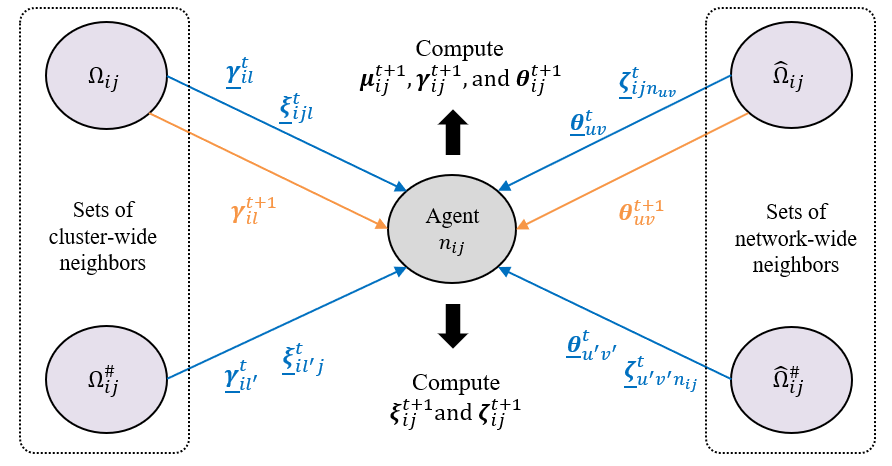}\\
  \caption{An illustrative data acquisition process of agent $n_{ij}$ (relabeled index). The arrows in {blue} represent the data flow by (\ref{2})-(\ref{3}) and the arrows in {orange} represent the data flow by (\ref{4})-(\ref{5}), $l \in \Omega_{ij},l^{'} \in \Omega^{\sharp}_{ij},n_{uv} \in \widehat{\Omega}_{ij},n_{u'v'} \in \widehat{\Omega}^{\sharp}_{ij}$. }\label{ccp}
\end{figure}
\end{Remark}

{{\begin{Remark}
Different from \cite{wang2023composite,wang2021distributed,wang2022composite}, the Asyn-DDPG algorithm can solve cluster-based affine-constrained optimization problems in a general incomplete graph with communication delays. Specifically, the algorithm in \cite{wang2023composite} can only be applied in a complete graph (or incomplete graph but no general affine constraints), and the analysis of the delayed information requires the common knowledge of time slots. In \cite{wang2022composite}, the communication graph should be complete for general affine constraints (or incomplete graph but no general affine constraints) and the length of delays is globally uniform. No communication delays are considered in \cite{wang2021distributed}. Moreover, all the algorithms in \cite{wang2023composite,wang2021distributed,wang2022composite} do not consider the cluster-based problem settings as this work.
\end{Remark}}}

{To further illustrate the notations involved in the Asyn-DDPG algorithm, the explanation of some important sets and variables is provided in Table \ref{tm2}.} 

{\begin{table}
{\caption{Meaning of Notations}\label{tm2}
\label{tab2}
\begin{center}
\begin{tabular}{p{3mm}|p{25mm}|p{90mm}}
\bottomrule
& Notation &  Explanation \\
\hline
1 & $\Omega_{ij}$ ($\Omega^{\sharp}_{ij}$) & Neighbor set of the $j$th agent in cluster $i$, whose indices are greater (smaller) than $j$. \\
 \hline
2 & $\widehat{\Omega}_{ij}$ ($\widehat{\Omega}^{\sharp}_{ij}$) & Neighbor set of the $n_{ij}$th (relabeled index) agent, whose indices are greater (smaller) than $n_{ij}$. \\
 \hline
3  & ${\bm{\mu}}_{ij}$, ${\bm{\gamma}}_{ij}$, ${\bm{\theta}}_{ij}$  & Variables maintained by the $j$th agent in cluster $i$. \\
 \hline
4 & ${\bm{\xi}}_{ijl}$ (${\bm{\zeta}}_{ijk}$) & Variable maintained by the $j$th agent in cluster $i$, where $l$ ($k$) is the index of its neighbor in $\Omega_{ij}$ ($\widehat{\Omega}_{ij}$). \\
\hline
4 & ${\bm{\xi}}_{ij}$ (${\bm{\zeta}}_{ij}$) & Variable maintained by the $j$th agent in cluster $i$, which is composed of $|\Omega_{ij}|$ ($|\widehat{\Omega}_{ij}|$) elements. \\
 \hline
5 & ${\bm{\xi}}_{i}$ (${\bm{\zeta}}_{i}$) & Variable maintained by the agents in cluster $i$. \\
  \hline
6 & ${\bm{\xi}}$ (${\bm{\zeta}}$) & Overall collection of ${\bm{\xi}}_{i}$ (${\bm{\zeta}}_i$). \\
  \hline
7 & $\bm{\alpha}_{ij}$ & Variable maintained by the $j$th agent in cluster $i$. \\
   \hline
8 & $\bm{\alpha}_{i}$  & Variable maintained by the agents in cluster $i$. \\
   \hline
9 & $\bm{\alpha}$  & Overall collection of ${\bm{\alpha}}_{i}$. \\
   \hline
10 & $\bm{\omega}$  & Overall collection of ${\bm{\xi}}$ and ${\bm{\zeta}}$. \\
\toprule
\end{tabular}
\end{center}}
\end{table}}

\section{Convergence Analysis}\label{sa5}

The convergence analysis of the proposed Asyn-DDPG algorithm is conducted in this section.

\begin{Lemma}\label{lam2}
By Assumption \ref{a1}, $\nabla_{\bm{\alpha}_{ij}}  s_{ij}(\bm{\alpha}_{ij})$ is Lipschitz continuous with constant $h_{ij} \triangleq \frac{ \| \mathbf{W}_{ij} \|_2^2 }{\sigma_{ij}}$, $\forall i \in V$, $j \in V_i$.
\end{Lemma}

See the proof in Appendix \ref{lam2p}.


Define $ \mathbf{C} \triangleq \mathrm{diag}[\mathbf{C}_1,\mathbf{C}_2,...,\mathbf{C}_{N}]$, $\mathbf{C}_i \triangleq \mathrm{diag}[c_{i1},c_{i2},...,c_{in_i}] \otimes \mathbf{I}_{M+n_iM+B}$, $\mathbf{S} \triangleq \mathrm{diag}[\mathbf{S}_1,...,$ $\mathbf{S}_{N},\widetilde{\mathbf{S}}_1,...,\widetilde{\mathbf{S}}_{N}]$, $ \mathbf{S}_i \triangleq  \mathrm{diag}[\mathbf{S}_{i1},\mathbf{S}_{i2},...,\mathbf{S}_{in_i}]$, $ \quad \mathbf{S}_{ij} \triangleq \pi_{ij}\mathbf{I}_{|\Omega_{ij}|n_iM}$, $ \widetilde{\mathbf{S}}_i \triangleq \mathrm{diag}[\widetilde{\mathbf{S}}_{i1},\widetilde{\mathbf{S}}_{i2},...,\widetilde{\mathbf{S}}_{in_i}]$, $\widetilde{\mathbf{S}}_{ij} \triangleq \pi_{ij}\mathbf{I}_{|\widehat{\Omega}_{ij}|B}$, $\mathbf{H} \triangleq \mathrm{diag}[\mathbf{H}_1,\mathbf{H}_2,...,\mathbf{H}_{N}]$, $\mathbf{H}_i \triangleq \mathrm{diag}[h_{i1},h_{i2},...,h_{in_i}] \otimes \mathbf{I}_{M+n_iM+B}$, $\mathbf{B} \triangleq (1+d)^2\mathbf{S}$, $\mathbf{Q} \triangleq (1+d)\mathbf{S}^{-1}$, $\Gamma \triangleq \sum_{i \in V} \sum_{j \in V_i} 4 \pi_{ij} d \iota^2$. Then we have the following two lemmas.



\begin{Lemma}\label{ath1}
By Algorithm 1, for certain step index $T \geq 2d +1 $, we have $
\sum_{t=0}^{T} (\| \bm{\alpha} - \bm{\alpha}^{t+1} \|^2_{\mathbf{Z}^{\top}\mathbf{S}\mathbf{Z}} - \| \bm{\alpha} -\underline{\bm{\alpha}}^t\|^2_{\mathbf{Z}^{\top}\mathbf{S}\mathbf{Z}}) \leq  \| \bm{\alpha} - \bm{\alpha}^{T+1} \|^2_{\mathbf{Z}^{\top}\mathbf{S}\mathbf{Z}} + \Gamma$, $\sum_{t=0}^{T}  \| \bm{\alpha}^{t+1} - \underline{\bm{\alpha}}^{t} \|^2_{ \mathbf{Z}^{\top}\mathbf{S}\mathbf{Z}}
\leq  \sum_{t=0}^{T} \| \bm{\alpha}^{t+1} - \bm{\alpha}^{t} \|^2_{\mathbf{Z}^{\top}\mathbf{B}\mathbf{Z}}$, and $\sum_{t=0}^{T}  (\| \bm{\omega} - \underline{\bm{\omega}}^{t} \|^2_{\mathbf{S}^{-1}} - \|\bm{\omega} - {\bm{\omega}}^{t+1} \|^2_{\mathbf{S}^{-1}} ) \leq \| \bm{\omega} - \bm{\omega}^{0} \|^2_{\mathbf{Q}}$.
\end{Lemma}

See the proof in Appendix \ref{ath1p}.

\begin{Lemma}\label{lass}
Let
\begin{align}\label{ss}
 \tau_{\mathrm{min}}(\mathbf{C}^{-1}) \geq  \tau_{\mathrm{max}}(\mathbf{H}) + 2\tau_{\mathrm{max}}(\mathbf{Z}^{\top}\mathbf{B}\mathbf{Z}).
\end{align}
Then $ \mathbf{C}^{-1}-\mathbf{H} - \mathbf{Z}^{\top}\mathbf{B}\mathbf{Z} \succeq 0$ and $\mathbf{C}^{-1} - \mathbf{Z}^{\top}\mathbf{S}\mathbf{Z} \succ 0$.
\end{Lemma}
See the proof in Appendix \ref{lassp}.

{To characterize the performance of the Asyn-DDPG algorithm, we formulate the Lagrangian function of Problem (P5) as $\mathcal{L}_H ({\bm{\alpha}},\bm{\omega}) \triangleq H ( {\bm{\alpha}}) + \bm{\omega}^{\top}\mathbf{Z}\bm{\alpha}$.} Define ${K}$ as the set of the {saddle points} of $\mathcal{L}_H$. After solving the saddle point of $\mathcal{L}_H$, we continue to investigate the saddle point property of $\mathcal{L}_F$: $\mathcal{L}_F(\mathbf{y}^*, \mathbf{e}^*, \bm{\mu}^*, \bm{\nu}^*,\bm{\phi}^*) \leq \mathcal{L}_F(\mathbf{y}, \mathbf{e}, \bm{\mu}^*, \bm{\nu}^*,\bm{\phi}^*)$, which gives $\mathbf{y}^* = \arg \min_{\mathbf{y}} \mathcal{L}_F(\mathbf{y}, \mathbf{e}, \bm{\mu}^*, \bm{\nu}^*,\bm{\phi}^*)$. Then given certain point $(\bm{\alpha}^*,\bm{\omega}^*) \in K$, we further have $\mathbf{y}_{ij}^* = \arg \min_{\mathbf{y}_{ij}} f_{ij}(\mathbf{y}_{ij}) + (\bm{\mu}_{ij}^{*} + \mathsf{L}^{i\top}_{j} \bm{\nu}^{*}_{i} + \mathsf{A}^{\top}_{ij} \bm{\phi}^* )^{\top} \mathbf{y}_{ij}
    =  \arg \min_{\mathbf{y}_{ij}} f_{ij}(\mathbf{y}_{ij}) + (\bm{\mu}_{ij}^{*} + \mathsf{L}^{i\top}_{j} \bm{\gamma}_{ij}^{*} + \mathsf{A}^{\top}_{ij} \bm{\theta}_{ij}^{*})^{\top} \mathbf{y}_{ij}$,
    where we use $\bm{\gamma}^*_{ij}= \bm{\nu}_i^{*}$ and $\bm{\theta}^*_{ij} = \bm{\phi}^*$, since $\bm{\gamma}^*_{ij}$ and $\bm{\theta}^*_{ij}$ are the consensual estimates of $\bm{\nu}_i^{*}$ and $\bm{\phi}^*$, respectively. Then the optimal solution to Problem (P1) can be $
    \mathbf{x}^*_i = \mathbf{y}^*_{i1} = ... = \mathbf{y}^*_{in_i}$, $\forall i \in V$.    In the following, we discuss how the solutions generated by (10) to (14) converge to the point in set $K$.
%
\begin{Theorem}\label{th2}
Suppose that Assumptions \ref{a0}-\ref{as5} hold and the step sizes are selected based on (\ref{ss}). By Algorithm 1, we have
\begin{align}\label{}
\mathbf{0}  \in & \lim_{t \rightarrow +\infty} (\partial_{\bm{\alpha}} H_r (\bm{\alpha}^{t+1}) + \nabla_{\bm{\alpha}} H_s (\bm{\alpha}^t) + \mathbf{Z}^{\top} \underline{\bm{\omega}}^t), \\
\mathbf{0} = &  \lim_{t \rightarrow +\infty} \mathbf{Z} \bm{\alpha}^{t+1}.
\end{align}
In addition, for any $(\bm{\alpha}^*,\bm{\omega}^*) \in {K}$ and $T \geq 2d +1 $, the convergence rate can be given by $|H(\overline{\bm{\alpha}}^{T+1}) - H(\bm{\alpha}^*)| \leq \mathcal{O}\left(\frac{\Xi}{T+1}\right), \| \bm{\omega}^* \| \|\mathbf{Z} \overline{\bm{\alpha}}^{T+1} \| \leq \mathcal{O} \left(\frac{\Xi}{T+1} \right)$, where $\overline{\bm{\alpha}}^{T+1} \triangleq \frac{1}{T+1}\sum_{t=0}^T {\bm{\alpha}}^{t+1}$ and $\Xi \triangleq \|  \bm{\omega}^* \|^2_{4\mathbf{Q}} + \| \bm{\omega}^0 \|^2_{\mathbf{Q}} + \|  \bm{\alpha}^* -\bm{\alpha}^0 \|^2_{\frac{1}{2}\mathbf{C}^{-1}} + \frac{\Gamma}{2}$.
\end{Theorem}

See the proof in Appendix \ref{th2p}. Note that to ensure (\ref{ss}) to hold, one can simply choose the step size $c_{ij} \in \left(0, \frac{1}{h + 2\tau_{\mathrm{max}}(\mathbf{Z}^{\top}\mathbf{B}\mathbf{Z})} \right]$, where $h = \max_{i \in V, j \in V_i} h_{ij}$.
{
\begin{Remark}
The problem formulation of this work fills the gap of the existing works in considering the cluster-based consensus in affine-constrained optimization problems. To solve this problem, the concept of partial and mixed consensus protocol is newly proposed, which can be more general than the conventional global consensus protocol with significant differences in terms of both the algorithm development and implementation. As a consequence, the Asyn-DDPG algorithm realizes a parallel but asynchronous computation manner for the agents by tolerating the communication delays in the network.
\end{Remark}}

\section{Numerical Simulation}\label{sa6}


In this simulation, we consider a social welfare optimization problem in a commodity market. In this market, we aim to supply certain amount of commodities to multiple consumer regions such that the utility function of the whole consumer community is optimized (transport cost is assumed to be negligible).
The utility function of different regions can be obtained by some learning machines based on the regional information, which can be settled by {{ensemble method}} \cite{guo2017learning}. However, due to the possibly large-scale data sets and privacy perseveration issues, it can be inefficient or even infeasible to transmit the data sets among the machines. Therefore, distributed optimization framework can be applied.

\subsection{Simulation Setup} 
\begin{figure}[htbp]
  \centering
  \includegraphics[width=10cm]{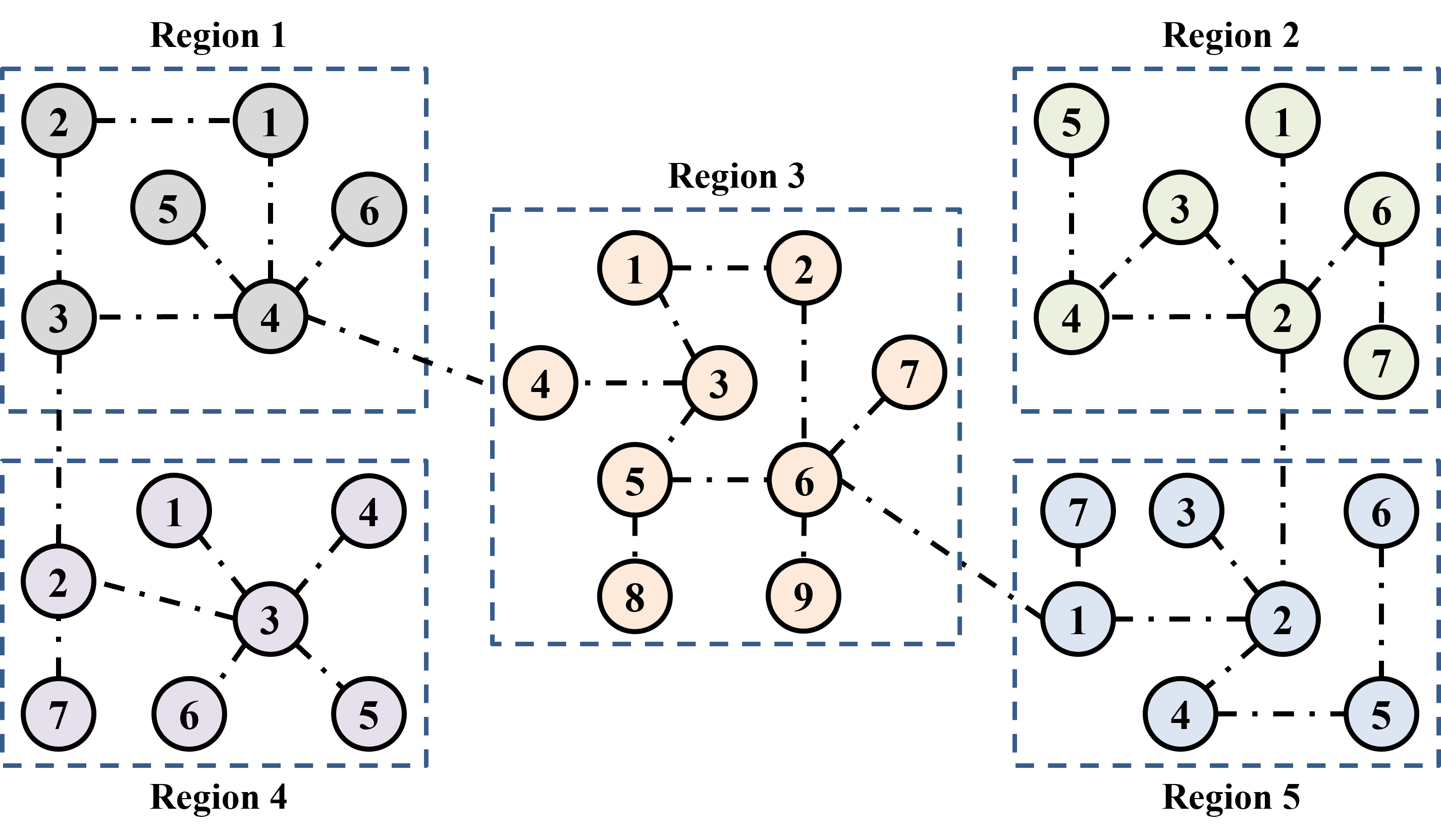}\\
  \caption{Communication typology of machines.}\label{mar2}
\end{figure}

Based on the discussion above, we let $f_{ij}(x_i)$ be the utility function generated by the $j$th machine in region $i$, where $x_i$ is the quantity of commodities supplied to region $i$. Then the utility function of region $i$ is settled by {{ensemble method}}, which is $f_i(x_i) \triangleq \frac{1}{n_i}\sum_{j\in V_i}f_{ij}(x_i)$ with $n_i$ being the number of machines in region $i$ \cite{guo2017learning}. Hence, the social welfare optimization problem can be formulated as
\begin{align}
  \mathrm{(P6)} \quad  \max\limits_{\mathbf{x} \in X} \quad  &  \sum_{i \in {V}} f_i (x_{i}) \quad
   \hbox{subject to} \quad  \mathbf{A}\mathbf{x} \leq b. \nonumber
\end{align}
Here, $\mathbf{x} \triangleq [x_1,...,x_{|V|}]^{\top}$, $\mathbf{A} \triangleq \mathbf{1}_{|V|}^{\top}$, and $b$ is total quantity of commodities in store. $X \triangleq \prod_{i\in V} \bigcap_{j \in V_i} X_{ij}$ with $X_{ij} \triangleq [\underline{x}_{ij},\overline{x}_{ij}]$ being the local constraint of the $j$th machine in region $i$. Then Problem (P6) follows the structure of Problem (P1) based on Remark \ref{r0}. In addition, the utility function obtained by the machines is assumed to be $f_{ij}(x_i) \triangleq -\varpi_{ij} x^2_i + \varsigma_{ij} x_i$ \cite{craven2005optimization}, where $\varpi_{ij}$ and $\varsigma_{ij}$ are set in Table \ref{tm31}. Let $\overline{x}_{ij} = \frac{\varsigma_{ij}}{2\varpi_{ij}}$, $\underline{x}_{ij} = 0$, $\forall i \in V, j \in V_i$, and $b=5$. The communication graph of the machines is shown in Fig. \ref{mar2}. The delay of communication links is selected from $\{0,1,2,...,10\}$ arbitrarily. To characterize the relative convergence error, we define $\varepsilon \triangleq \left|\frac{ H(\bm{\alpha}) - H(\bm{\alpha}^*) }{H(\bm{\alpha}^*)}\right|$ with a non-zero $H(\bm{\alpha}^*)$. {Note that the existing cluster-based algorithms cannot be directly applied in this simulation. For example, in \cite{guo2017distributed,li2019gossip,shi2019multi}, some special communication graph topologies are required and no communication delays are considered. In addition, other existing algorithms as well as their variants still cannot be applied due to the different problem settings. See detailed discussions in Section I-B.}

\begin{table}[H]
\caption{Parameters of Simulation}\label{tm31}
\label{tab2}
\begin{center}
\begin{tabular}{p{4.5mm}|p{4.5mm}|p{4.5mm}|p{4.5mm}|p{4.5mm}|p{4.5mm}|p{4.5mm}|p{4.5mm}|p{4.5mm}|p{4.5mm}}
\bottomrule
$j$ & 1 & 2 & 3 & 4 & 5 & 6 & 7 & 8 & 9 \\
\hline
$\varpi_{1j}$ & 0.8 & 0.6 & 0.7 & 0.9 & 1.0 & 0.7 & - & - & - \\
\hline
$\varsigma_{1j}$ & 3.5 & 3.2 & 3.1 & 3.4 & 3.6 & 3.3 & - & - & - \\
\hline
$\varpi_{2j}$ & 1.4 & 1.5 & 1.2 & 1.6 & 1.8 & 1.9 & 1.5 & - & - \\
\hline
$\varsigma_{2j}$ & 2.7 & 2.8 & 2.6 & 2.8 & 2.9 & 3.0 & 2.7 & - & - \\
\hline
$\varpi_{3j}$ & 0.9 & 1.2 & 1.5 & 1.6 & 0.8 & 1.0 & 2.2 & 2.5 & 1.7 \\
\hline
$\varsigma_{3j}$ & 4.1 & 4.2 & 4.3 & 4.2 & 4.5 & 4.4 & 4.7 & 4.1 & 4.2 \\
\hline
$\varpi_{4j}$ & 0.5 & 0.6 & 0.7 & 0.6 & 0.5 & 0.8 & 0.9 & - & - \\
\hline
$\varsigma_{4j}$ & 1.7 & 1.8 & 1.6 & 1.9 & 2.0 & 1.8 & 1.5 & - & - \\
\hline
$\varpi_{5j}$ & 1.2 & 1.1 & 1.3 & 0.9 & 1.0 & 1.2 & 1.1 & - & - \\
\hline
$\varsigma_{5j}$ & 2.2 & 2.3 & 1.4 & 2.1 & 1.9 & 2.5 & 2.4 & - & - \\
\toprule
\end{tabular}
\end{center}
\end{table}

\subsection{Simulation Result} 
\begin{figure*}
\centering
\subfigure[Trajectory of $\bm{\alpha}$.]{
\begin{minipage}[t]{0.33\linewidth}
\centering
\includegraphics[width=8cm,height=6cm]{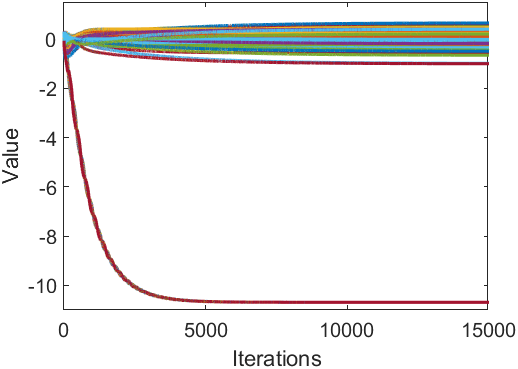}
\end{minipage}%
} \\%
\subfigure[Trajectory of $\bm{\omega}$.]{
\begin{minipage}[t]{0.33\linewidth}
\centering
\includegraphics[width=8cm,height=6cm]{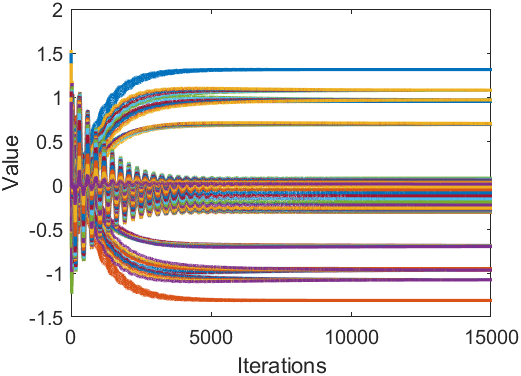}
\end{minipage}%
} \\%
\subfigure[Trajectory of $\varepsilon$.]{
\begin{minipage}[t]{0.33\linewidth}
\centering
\includegraphics[width=8cm,height=6cm]{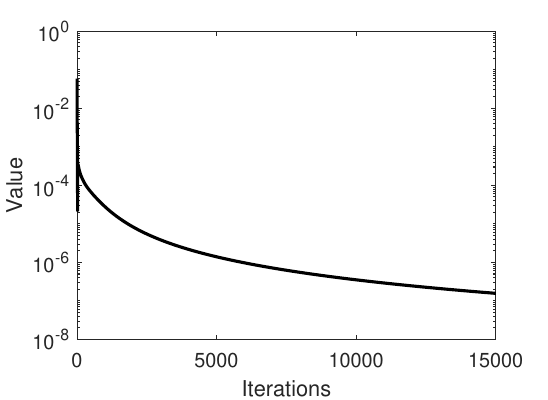}
\end{minipage}%
}%
\caption{Simulation result.}\label{g3}
\end{figure*}

The simulation result is shown in Fig. \ref{g3}. Fig. \ref{g3}-(a) shows that $\bm{\alpha}$ tends to the steady state gradually. Fig. \ref{g3}-(b) depicts the trajectory of $\bm{\omega}$. The trajectory of convergence error is shown in Fig. \ref{g3}-(c). Finally, with the achieved optimal $\bm{\alpha}$ and $\bm{\omega}$, the optimal solution to Problem (P6) is obtained as $\mathbf{x}^*= [0.38, 0.89, 1.44, 1.34, 0.95]^{\top}$.

\section{Conclusion}\label{sa7}

In this work, we considered a novel DOP in a multi-cluster network with affine coupling constraints. In this problem, each cluster can make its decision based on the consensus protocol among the agents involved. To achieve the optimal solution of the whole network, an Asyn-DDPG algorithm was proposed, where the agents are only required to access the information from their neighbors with communication delays. The example in the simulation demonstrates the effectiveness of the Asyn-DDPG algorithm in solving some practical problems in the industry.

\section{Appendix}

\subsection{Proof of Lemma \ref{lam2}}\label{lam2p}

By \cite[Lemma V.7]{notarnicola2016asynchronous}, $ \nabla f^{\circ}_{ij}$ is Lipschitz continuous with Lipschitz constant $\frac{1}{\sigma_{ij}}$. Then by \cite[Lemma 3.1]{beck2014fast}, $ \nabla_{\bm{\alpha}_{ij}} f^{\circ}_{ij}(\mathbf{W}_{ij}\bm{\alpha}_{ij})$ is Lipschitz continuous with constant $\frac{ \| \mathbf{W}_{ij} \|_2^2 }{\sigma_{ij}}$, which means $\nabla_{\bm{\alpha}_{ij}} s_{ij} (\bm{\alpha}_{ij})=\nabla_{\bm{\alpha}_{ij}} f^{\circ}_{ij}(\mathbf{W}_{ij}\bm{\alpha}_{ij}) + \mathbf{D}^{\top}_{ij}$ is Lipschitz continuous with constant $\frac{ \| \mathbf{W}_{ij} \|_2^2 }{\sigma_{ij}}$.

\subsection{{Proof of Lemma \ref{ath1}}}\label{ath1p}

To facilitate the derivation, define $\bm{\varphi}_{ij} \triangleq \mathbf{Z}_{ij}\bm{\alpha}$, $\bm{\varphi}^t_{ij} \triangleq \mathbf{Z}_{ij}\bm{\alpha}^t$, and $\underline{\bm{\varphi}}_{ij}^t \triangleq \mathbf{Z}_{ij}\underline{\bm{\alpha}}^t = \mathbf{Z}_{ij}{\bm{\alpha}}^{(t-d)^+} = {\bm{\varphi}}_{ij}^{(t-d)^+}$, $\forall i \in V$, $j\in V_i$. Therefore, $\| \bm{\varphi}_{ij}\|\leq \iota$. Then we can have the following results.

\begin{align}\label{54}
&\sum_{t=0}^{T} (\| \bm{\alpha} - \bm{\alpha}^{t+1} \|^2_{\mathbf{Z}^{\top}\mathbf{S}\mathbf{Z}} - \| \bm{\alpha} -\underline{\bm{\alpha}}^t\|^2_{\mathbf{Z}^{\top}\mathbf{S}\mathbf{Z}}) \nonumber \\
= & \sum_{t=0}^{T}  (\| \mathbf{Z}(\bm{\alpha} - {\bm{\alpha}}^{t+1}) \|^2_{\mathbf{S}} - \| \mathbf{Z}(\bm{\alpha} - \underline{\bm{\alpha}}^{t}) \|^2_{\mathbf{S}})\nonumber \\
= & \sum_{i \in V} \sum_{j \in V_i} \pi_{ij} \sum_{t=0}^{T}  (\| \mathbf{Z}_{ij}(\bm{\alpha} - {\bm{\alpha}}^{t+1}) \|^2 - \| \mathbf{Z}_{ij}(\bm{\alpha} - \underline{\bm{\alpha}}^{t}) \|^2 )\nonumber \\
 = &\sum_{i \in V} \sum_{j \in V_i} \pi_{ij} \sum_{t=0}^{T} ( \| \bm{\varphi}_{ij} - {\bm{\varphi}}_{ij}^{t+1} \|^2 - \| \bm{\varphi}_{ij} - \underline{\bm{\varphi}}_{ij}^{t} \|^2 ) \nonumber \\
\stackrel{(i)}{=} &  \sum_{i \in V} \sum_{j \in V_i} \pi_{ij} (\| \bm{\varphi}_{ij} - \bm{\varphi}_{ij}^{1} \|^2- \| \bm{\varphi}_{ij} - \bm{\varphi}_{ij}^{0} \|^2 \nonumber \\
 & + \| \bm{\varphi}_{ij} - \bm{\varphi}_{ij}^{2} \|^2 - \| \bm{\varphi}_{ij} - \bm{\varphi}_{ij}^{0}\|^2  + \cdots \nonumber \\
 & + \| \bm{\varphi}_{ij} - \bm{\varphi}_{ij}^{d+1} \|^2 - \| \bm{\varphi}_{ij} - \bm{\varphi}_{ij}^{0} \|^2 \nonumber \\
 & + \| \bm{\varphi}_{ij} - \bm{\varphi}_{ij}^{d+2} \|^2- \| \bm{\varphi}_{ij} - \bm{\varphi}_{ij}^{1}\|^2 \nonumber \\
 & + \| \bm{\varphi}_{ij} - \bm{\varphi}_{ij}^{d+3} \|^2 - \| \bm{\varphi}_{ij} - \bm{\varphi}_{ij}^{2} \|^2 + \cdots \nonumber \\
  &  + \| \bm{\varphi}_{ij} - \bm{\varphi}_{ij}^{T-d} \|^2 - \| \bm{\varphi}_{ij} - \bm{\varphi}_{ij}^{T-2d-1}\|^2 \nonumber \\
  &  + \| \bm{\varphi}_{ij} - \bm{\varphi}_{ij}^{T-d+1} \|^2 - \| \bm{\varphi}_{ij} - \bm{\varphi}_{ij}^{T-2d}\|^2 + \cdots \nonumber \\
 &  + \| \bm{\varphi}_{ij} - \bm{\varphi}_{ij}^{T} \|^2 - \| \bm{\varphi}_{ij} - \bm{\varphi}_{ij}^{T-d-1}\|^2  \nonumber \\
 &  + \| \bm{\varphi}_{ij} - \bm{\varphi}_{ij}^{T+1} \|^2 - \| \bm{\varphi}_{ij} - \bm{\varphi}_{ij}^{T-d}\|^2 ) \nonumber \\
\stackrel{(ii)}{=} & \sum_{i \in V} \sum_{j \in V_i} \pi_{ij} \bigg(\sum_{t=T-d+1}^{T+1} \| \bm{\varphi}_{ij} - \bm{\varphi}_{ij}^{t} \|^2 \nonumber \\
 & \quad \quad \quad \quad \quad \quad -(d+1) \| \bm{\varphi}_{ij} - \bm{\varphi}_{ij}^{0} \|^2 \bigg)  \nonumber \\
  \leq & \sum_{i \in V} \sum_{j \in V_i} \pi_{ij} \bigg(\sum_{t=T-d+1}^{T} (2\| \bm{\varphi}_{ij} \|^2 + 2\|\bm{\varphi}_{ij}^{t} \|^2)  \nonumber \\
 & \quad \quad \quad \quad + \| \bm{\varphi}_{ij} - \bm{\varphi}_{ij}^{T+1} \|^2 -(d+1) \| \bm{\varphi}_{ij} - \bm{\varphi}_{ij}^{0} \|^2 \bigg)  \nonumber \\
\leq & \sum_{i \in V} \sum_{j \in V_i} \pi_{ij} (4d \iota^2 + \| \bm{\varphi}_{ij} - \bm{\varphi}_{ij}^{T+1} \|^2)  \nonumber \\
 = & \| \bm{\alpha} - \bm{\alpha}^{T+1} \|^2_{\mathbf{Z}^{\top}\mathbf{S}\mathbf{Z}} +\sum_{i \in V} \sum_{j \in V_i} 4\pi_{ij} d \iota^2,
\end{align}
where $(i)$ holds by expanding the summation from $t=0$ to $T$ and $(ii)$ holds by performing cancellations.

\begin{align}\label{55}
& \sum_{t=0}^{T}  \| \bm{\alpha}^{t+1} - \underline{\bm{\alpha}}^{t} \|^2_{ \mathbf{Z}^{\top}\mathbf{S}\mathbf{Z}} \nonumber \\
= & \sum_{t=0}^{T}  \| \mathbf{Z}(\bm{\alpha}^{t+1} - \underline{\bm{\alpha}}^{t}) \|^2_{ \mathbf{S}} \nonumber \\
 = &\sum_{i \in V} \sum_{j \in V_i} \pi_{ij} \sum_{t=0}^{T} \| \mathbf{Z}_{ij} (\bm{\alpha}^{t+1} - \underline{\bm{\alpha}}^{t}) \|^2  \nonumber \\
 = &\sum_{i \in V} \sum_{j \in V_i} \pi_{ij} \sum_{t=0}^{T} \| \bm{\varphi}_{ij}^{t+1} - \underline{\bm{\varphi}}_{ij}^{t} \|^2  \nonumber \\
\stackrel{(i)}{\leq} & \sum_{i \in V} \sum_{j \in V_i} \pi_{ij} \sum_{t=0}^{T}  (d+1) (\| \bm{\varphi}_{ij}^{t+1} - \bm{\varphi}_{ij}^{t} \|^2  \nonumber \\
& + \| \bm{\varphi}_{ij}^{t} - \bm{\varphi}_{ij}^{t-1} \|^2 + \cdots + \| \bm{\varphi}_{ij}^{(t-d)^++1} - \bm{\varphi}_{ij}^{(t-d)^+} \|^2 ) \nonumber \\
\stackrel{(ii)}{=} & \sum_{i \in V} \sum_{j \in V_i} \pi_{ij} (d+1) ((\| \bm{\varphi}_{ij}^{T+1} - \bm{\varphi}_{ij}^{T} \|^2  \nonumber \\
& + \| \bm{\varphi}_{ij}^{T} - \bm{\varphi}_{ij}^{T-1} \|^2 + \cdots + \| \bm{\varphi}_{ij}^{T-d+1} - \bm{\varphi}_{ij}^{T-d} \|^2) \nonumber \\
& + (\| \bm{\varphi}_{ij}^{T} - \bm{\varphi}_{ij}^{T-1} \|^2 + \cdots +  \| \bm{\varphi}_{ij}^{T-d} - \bm{\varphi}_{ij}^{T-d-1} \|^2 ) \nonumber \\
& + \cdots + ( \| \bm{\varphi}_{ij}^{d+1} - \bm{\varphi}_{ij}^{d} \|^2  + \cdots + \| \bm{\varphi}_{ij}^{1} - \bm{\varphi}_{ij}^{0} \|^2 ) + \cdots \nonumber \\
& + ( \| \bm{\varphi}_{ij}^{2} - \bm{\varphi}_{ij}^{1} \|^2 + \| \bm{\varphi}_{ij}^{1} - \bm{\varphi}_{ij}^{0} \|^2 )  + \| \bm{\varphi}_{ij}^{1} - \bm{\varphi}_{ij}^{0} \|^2 ) \nonumber \\
= & \sum_{i \in V} \sum_{j \in V_i} \pi_{ij} (d+1) (\| \bm{\varphi}_{ij}^{T+1} - \bm{\varphi}_{ij}^{T} \|^2 \nonumber \\
& + 2 \| \bm{\varphi}_{ij}^{T} - \bm{\varphi}_{ij}^{T-1} \|^2 + 3 \| \bm{\varphi}_{ij}^{T-1} - \bm{\varphi}_{ij}^{T-2} \|^2 + \cdots \nonumber \\
& + (d+1) \| \bm{\varphi}_{ij}^{T-d+1} - \bm{\varphi}_{ij}^{T-d} \|^2 \nonumber \\
& + (d+1) \| \bm{\varphi}_{ij}^{T-d} - \bm{\varphi}_{ij}^{T-d-1} \|^2  \nonumber \\
& + (d+1) \| \bm{\varphi}_{ij}^{T-d-1} - \bm{\varphi}_{ij}^{T -d-2} \|^2 + \cdots \nonumber \\
& + (d+1) \| \bm{\varphi}_{ij}^{2} - \bm{\varphi}_{ij}^{1} \|^2  + (d+1) \| \bm{\varphi}_{ij}^{1} - \bm{\varphi}_{ij}^{0} \|^2 ) \nonumber \\
\leq &  \sum_{t=0}^{T} \sum_{i \in V} \sum_{j \in V_i} \pi_{ij}  (d+1)^2 \| \bm{\varphi}_{ij}^{t+1} - \bm{\varphi}_{ij}^{t} \|^2 \nonumber \\
= & \sum_{t=0}^{T} \| \bm{\alpha}^{t+1} - \bm{\alpha}^{t} \|^2_{\mathbf{Z}^{\top}\mathbf{B}\mathbf{Z}},
\end{align}
where $(i)$ is based on {\em{Cauchy-Schwarz inequality}} and $(ii)$ holds by expanding the summation from $t=0$ to $T$.

\begin{align}\label{}
& \sum_{t=0}^{T}  ( \|\bm{\omega} - \underline{\bm{\omega}}^{t} \|^2_{\mathbf{S}^{-1}} - \| \bm{\omega} - {\bm{\omega}}^{t+1} \|^2_{\mathbf{S}^{-1}} )\nonumber \\
 = &\sum_{i \in V} \sum_{j \in V_i} \frac{1}{\pi_{ij}} \sum_{t=0}^{T} ( \| \bm{\xi}_{ij} - \underline{\bm{\xi}}_{ij}^{t} \|^2 - \| \bm{\xi}_{ij} - {\bm{\xi}}_{ij}^{t+1} \|^2 ) \nonumber \\
 & + \sum_{i \in V} \sum_{j \in V_i} \frac{1}{\pi_{ij}} \sum_{t=0}^{T} (  \| \bm{\zeta}_{ij} - \underline{\bm{\zeta}}_{ij}^{t} \|^2 - \| \bm{\zeta}_{ij} - {\bm{\zeta}}_{ij}^{t+1} \|^2) \nonumber \\
 \stackrel{(i)}{=} &  \sum_{i \in V} \sum_{j \in V_i} \frac{1}{\pi_{ij}} (\| \bm{\xi}_{ij} - \bm{\xi}_{ij}^{0} \|^2- \| \bm{\xi}_{ij} - \bm{\xi}_{ij}^{1} \|^2 \nonumber \\
 & + \| \bm{\xi}_{ij} - \bm{\xi}_{ij}^{0} \|^2 - \| \bm{\xi}_{ij} - \bm{\xi}_{ij}^{2}\|^2  + \cdots \nonumber \\
 & + \| \bm{\xi}_{ij} - \bm{\xi}_{ij}^{0} \|^2 - \| \bm{\xi}_{ij} - \bm{\xi}_{ij}^{d+1} \|^2 \nonumber \\
 & + \| \bm{\xi}_{ij} - \bm{\xi}_{ij}^{1} \|^2- \| \bm{\xi}_{ij} - \bm{\xi}_{ij}^{d+2}\|^2 \nonumber \\
 & + \| \bm{\xi}_{ij} - \bm{\xi}_{ij}^{2} \|^2 - \| \bm{\xi}_{ij} - \bm{\xi}_{ij}^{d+3} \|^2 + \cdots \nonumber \\
  &  + \| \bm{\xi}_{ij} - \bm{\xi}_{ij}^{T-2d-1} \|^2 - \| \bm{\xi}_{ij} - \bm{\xi}_{ij}^{T-d}\|^2 \nonumber \\
  &  + \| \bm{\xi}_{ij} - \bm{\xi}_{ij}^{T-2d} \|^2 - \| \bm{\xi}_{ij} - \bm{\xi}_{ij}^{T-d+1}\|^2 + \cdots \nonumber \\
 &  + \| \bm{\xi}_{ij} - \bm{\xi}_{ij}^{T-d-1} \|^2 - \| \bm{\xi}_{ij} - \bm{\xi}_{ij}^{T}\|^2  \nonumber \\
 &  + \| \bm{\xi}_{ij} - \bm{\xi}_{ij}^{T-d} \|^2 - \| \bm{\xi}_{ij} - \bm{\xi}_{ij}^{T+1}\|^2 ) \nonumber \\
 & + \sum_{i \in V} \sum_{j \in V_i} \frac{1}{\pi_{ij}} (\| \bm{\zeta}_{ij} - \bm{\zeta}_{ij}^{0} \|^2- \| \bm{\zeta}_{ij} - \bm{\zeta}_{ij}^{1} \|^2 \nonumber \\
 & + \| \bm{\zeta}_{ij} - \bm{\zeta}_{ij}^{0} \|^2 - \| \bm{\zeta}_{ij} - \bm{\zeta}_{ij}^{2}\|^2  + \cdots \nonumber \\
 & + \| \bm{\zeta}_{ij} - \bm{\zeta}_{ij}^{0} \|^2 - \| \bm{\zeta}_{ij} - \bm{\zeta}_{ij}^{d+1} \|^2 \nonumber \\
 & + \| \bm{\zeta}_{ij} - \bm{\zeta}_{ij}^{1} \|^2- \| \bm{\zeta}_{ij} - \bm{\zeta}_{ij}^{d+2}\|^2 \nonumber \\
 & + \| \bm{\zeta}_{ij} - \bm{\zeta}_{ij}^{2} \|^2 - \| \bm{\zeta}_{ij} - \bm{\zeta}_{ij}^{d+3} \|^2 + \cdots \nonumber \\
  &  + \| \bm{\zeta}_{ij} - \bm{\zeta}_{ij}^{T-2d-1} \|^2 - \| \bm{\zeta}_{ij} - \bm{\zeta}_{ij}^{T-d}\|^2 \nonumber \\
  &  + \| \bm{\zeta}_{ij} - \bm{\zeta}_{ij}^{T-2d} \|^2 - \| \bm{\zeta}_{ij} - \bm{\zeta}_{ij}^{T-d+1}\|^2 + \cdots \nonumber \\
 &  + \| \bm{\zeta}_{ij} - \bm{\zeta}_{ij}^{T-d-1} \|^2 - \| \bm{\zeta}_{ij} - \bm{\zeta}_{ij}^{T}\|^2  \nonumber \\
 &  + \| \bm{\zeta}_{ij} - \bm{\zeta}_{ij}^{T-d} \|^2 - \| \bm{\zeta}_{ij} - \bm{\zeta}_{ij}^{T+1}\|^2 ) \nonumber \\
\stackrel{(ii)}{=} & \sum_{i \in V} \sum_{j \in V_i} \frac{1}{\pi_{ij}} \bigg((d+1) \| \bm{\xi}_{ij} - \bm{\xi}_{ij}^{0} \|^2  -\sum_{t=T-d+1}^{T+1} \| \bm{\xi}_{ij} - \bm{\xi}_{ij}^{t} \|^2 \bigg)  \nonumber \\
 & + \sum_{i \in V} \sum_{j \in V_i} \frac{1}{\pi_{ij}} \bigg((d+1) \| \bm{\zeta}_{ij} - \bm{\zeta}_{ij}^{0} \|^2  -\sum_{t=T-d+1}^{T+1} \| \bm{\zeta}_{ij} - \bm{\zeta}_{ij}^{t} \|^2 \bigg)  \nonumber \\
\leq &  \sum_{i \in V} \sum_{j \in V_i} \frac{d+1}{\pi_{ij}}  (\| \bm{\xi}_{ij} - \bm{\xi}_{ij}^{0} \|^2  +  \| \bm{\zeta}_{ij} - \bm{\zeta}_{ij}^{0} \|^2) \nonumber \\
= & \| \bm{\omega} - \bm{\omega}^{0} \|^2_{\mathbf{Q}},
\end{align}
where $(i)$ and $(ii)$ are obtained by expanding the summation according to $t$ and performing cancellations, respectively.

\subsection{Proof of Lemma \ref{lass}}\label{lassp}
Note that $\mathbf{B} - \mathbf{S}$ is a diagonal matrix with elements $(d^2+2d)\pi_{ij}$ on the diagonal, $i\in V$, $j \in V_i$. Therefore, $\mathbf{B} - \mathbf{S} \succeq 0$. Then by \cite[Thm. 3.1]{horn1998eigenvalue}, we have
\begin{align}\label{}
& \tau_{\mathrm{max}}(\mathbf{Z}^{\top}\mathbf{B}\mathbf{Z}) - \tau_{\mathrm{max}}(\mathbf{Z}^{\top}\mathbf{S}\mathbf{Z}) \nonumber \\
= & \tau_{\mathrm{max}}(\mathbf{Z}^{\top}\mathbf{B}\mathbf{Z}) + \tau_{\mathrm{min}}(- \mathbf{Z}^{\top}\mathbf{S}\mathbf{Z}) \nonumber \\
\geq & \tau_{\mathrm{min}}(\mathbf{Z}^{\top}\mathbf{B}\mathbf{Z} - \mathbf{Z}^{\top}\mathbf{S}\mathbf{Z}) \nonumber \\
= & \tau_{\mathrm{min}} (\mathbf{Z}^{\top}(\mathbf{B}-\mathbf{S})\mathbf{Z}) \nonumber \\
= & \tau_{\mathrm{min}}(((\mathbf{B}-\mathbf{S})^{\frac{1}{2}}\mathbf{Z})^{\top} (\mathbf{B}-\mathbf{S})^{\frac{1}{2}}\mathbf{Z}) \geq 0.
\end{align}

On the other hand, note that $\mathbf{C}^{-1}$, $\mathbf{H}$, $\mathbf{Z}^{\top}\mathbf{B}\mathbf{Z}$, and $\mathbf{Z}^{\top}\mathbf{S}\mathbf{Z}$ are Hermitian. Based on \cite[Thm. 3.1]{horn1998eigenvalue} and (\ref{ss}), we have
\begin{align}\label{}
0 \leq & \tau_{\mathrm{min}}(\mathbf{C}^{-1}) - \tau_{\mathrm{max}}(\mathbf{H}) - 2 \tau_{\mathrm{max}}(\mathbf{Z}^{\top}\mathbf{B}\mathbf{Z}) \nonumber \\
\stackrel{(i)}{\leq} & \tau_{\mathrm{min}}(\mathbf{C}^{-1}) + \tau_{\mathrm{min}}(- \mathbf{H}) + \tau_{\mathrm{min}}( - \mathbf{Z}^{\top}\mathbf{B}\mathbf{Z}) \nonumber \\
\leq & \tau_{\mathrm{min}}(\mathbf{C}^{-1} - \mathbf{H} - \mathbf{Z}^{\top}\mathbf{B}\mathbf{Z})
\end{align}
and
\begin{align}\label{}
0 \leq & \tau_{\mathrm{min}}(\mathbf{C}^{-1}) - \tau_{\mathrm{max}}(\mathbf{H}) - 2 \tau_{\mathrm{max}}(\mathbf{Z}^{\top}\mathbf{B}\mathbf{Z}) \nonumber \\
\stackrel{(ii)}{<} & \tau_{\mathrm{min}}(\mathbf{C}^{-1}) - \tau_{\mathrm{max}}( \mathbf{Z}^{\top}\mathbf{S}\mathbf{Z}) \nonumber \\
= & \tau_{\mathrm{min}}(\mathbf{C}^{-1}) + \tau_{\mathrm{min}}(-\mathbf{Z}^{\top}\mathbf{S}\mathbf{Z}) \nonumber \\
\leq & \tau_{\mathrm{min}}(\mathbf{C}^{-1} - \mathbf{Z}^{\top}\mathbf{S}\mathbf{Z}),
\end{align}
where $(i)$ holds since $\mathbf{Z}^{\top}\mathbf{B}\mathbf{Z}\succeq 0$, $(ii)$ holds with $\tau_{\mathrm{max}}(\mathbf{Z}^{\top}\mathbf{B}\mathbf{Z}) - \tau_{\mathrm{max}}(\mathbf{Z}^{\top}\mathbf{S}\mathbf{Z})  \geq 0$ and $\mathbf{H} \succ 0$. Hence, $\mathbf{C}^{-1}-\mathbf{H} - \mathbf{Z}^{\top}\mathbf{B}\mathbf{Z} \succeq 0$ and $\mathbf{C}^{-1} - \mathbf{Z}^{\top}\mathbf{S}\mathbf{Z} \succ 0$.

\subsection{Proof of Theorem \ref{th2}}\label{th2p}

Since $s_{ij}$ is convex and Lipschitz continuously differentiable, we can have
\begin{align}\label{e3}
0 = & \sum_{i \in {V}} \sum_{j \in {V}_i} (\bm{\alpha}_{ij} - \bm{\alpha}_{ij}^t)^{\top} \nabla_{\bm{\alpha}_{ij}} s_{ij} (\bm{\alpha}_{ij}^t) \nonumber \\
& + \sum_{i \in {V}} \sum_{j \in {V}_i}(\bm{\alpha}_{ij}^t - \bm{\alpha}_{ij}^{t+1})^{\top} \nabla_{\bm{\alpha}_{ij}} s_{ij} (\bm{\alpha}_{ij}^t) \nonumber \\
& - \sum_{i \in {V}} \sum_{j \in {V}_i}(\bm{\alpha}_{ij} - \bm{\alpha}_{ij}^{t+1})^{\top} \nabla_{\bm{\alpha}_{ij}} s_{ij} (\bm{\alpha}_{ij}^t) \nonumber \\
\leq & \sum_{i \in {V}} \sum_{j \in {V}_i} ( s_{ij}(\bm{\alpha}_{ij}) -  s_{ij}(\bm{\alpha}_{ij}^t) )   \nonumber \\
&  +  \sum_{i \in {V}} \sum_{j \in {V}_i} ( s_{ij}(\bm{\alpha}_{ij}^t)-  s_{ij}(\bm{\alpha}_{ij}^{t+1}) ) \nonumber \\
& + \sum_{i \in {V}} \sum_{j \in {V}_i} \frac{h_{ij}}{2} \|  \bm{\alpha}_{ij}^t - \bm{\alpha}_{ij}^{t+1} \|^2 \nonumber  \\
& - \sum_{i \in {V}} \sum_{j \in {V}_i}(\bm{\alpha}_{ij} - \bm{\alpha}_{ij}^{t+1})^{\top} \nabla_{\bm{\alpha}_{ij}} s_{ij} (\bm{\alpha}_{ij}^t) \nonumber \\
= & H_s(\bm{\alpha})  -  H_s(\bm{\alpha}^{t+1})  + \| \bm{\alpha}^t  - \bm{\alpha}^{t+1} \|^2_{\frac{1}{2}\mathbf{H}} \nonumber  \\
& - (  \bm{\alpha} - \bm{\alpha}^{t+1})^{\top}  \nabla_{\bm{\alpha}} H_s (\bm{\alpha}^t).
\end{align}
Note that the proximal mapping of an indicator function is equivalent to a Euclidean projection \cite[Sec. 1.2]{parikh2014proximal}, which means $ \mathcal{P}_{Y_i}  = \mathbf{prox}^{c_{ij}}_{\mathbb{I}_{Y_i}}$ and $
 \mathcal{P}_{J} = \mathbf{prox}^{c_{ij}}_{\mathbb{I}_{J}}$. In addition, by $r_{ij} (\bm{\alpha}_{ij}) = g_{ij}^{\circ}(\bm{\mu}_{ij}) + \mathbb{I}_{Y_i}(\bm{\gamma}_{ij}) + \mathbb{I}_{J}(\bm{\theta}_{ij})$, we can have \cite[Thm. 6.6]{beck2017first}
\begin{align}\label{m7}
\mathbf{prox}^{c_{ij}}_{r_{ij}} = \mathbf{prox}^{c_{ij}}_{g_{ij}^{\circ}} \times \mathbf{prox}^{c_{ij}}_{\mathbb{I}_{Y_i}} \times \mathbf{prox}^{c_{ij}}_{\mathbb{I}_{J}}.
\end{align}

By \cite[Lemma V.7]{notarnicola2016asynchronous}, we have
\begin{align}\label{66d}
  \nabla_{\bm{\alpha}_{ij}} s_{ij} (\bm{\alpha}^t_{ij})
    = & \mathbf{W}_{ij}^{\top}  \nabla_{\mathbf{W}_{ij} \bm{\alpha}_{ij}} f^{\circ}_{ij}(\mathbf{W}_{ij} \bm{\alpha}^t_{ij}) + \mathbf{D}_{ij}^{\top} \nonumber \\
     = & \mathbf{W}_{ij}^{\top} \arg \max\limits_{\mathbf{n}} ( (\mathbf{W}_{ij} \bm{\alpha}^t_{ij})^{\top} \mathbf{n} - f_{ij}(\mathbf{n})) + \mathbf{D}_{ij}^{\top}.
\end{align}
Then based on (\ref{m1})-(\ref{m3}), it can be verified that $\nabla_{\bm{\alpha}_{ij}} s_{ij} (\bm{\alpha}^t_{ij}) = [\mathbf{m}_{\bm{\mu}_{ij}}^{t\top}, \mathbf{m}_{\bm{\gamma}_{ij}}^{t\top}, \mathbf{m}_{\bm{\theta}_{ij}}^{t\top}]^{\top}$, which means
\begin{align}
& \mathbf{m}_{\bm{\mu}_{ij}}^t = \nabla_{\bm{\mu}_{ij}} s_{ij} (\bm{\alpha}_{ij}^t), \label{m4} \\
& \mathbf{m}_{\bm{\gamma}_{ij}}^t = \nabla_{\bm{\gamma}_{ij}} s_{ij} (\bm{\alpha}_{ij}^t), \label{m5} \\
& \mathbf{m}_{\bm{\theta}_{ij}}^t = \nabla_{\bm{\theta}_{ij}} s_{ij} (\bm{\alpha}_{ij}^t).\label{m6}
\end{align}
Then (\ref{1}) to (\ref{3}) can be written into a compact form by including all $i \in V$ and $j \in V_i$ with the help of (\ref{m7}), which gives
\begin{align}
 \bm{\alpha}^{t+1} = & \mathbf{prox}^{\mathbf{C}}_{H_r}[\bm{\alpha}^t - \mathbf{C} ( \nabla_{\bm{\alpha}} H_s (\bm{\alpha}^t)  + {\mathbf{Z}}^{\top} \underline{\bm{\omega}}^t  + \mathbf{Z}^{\top}\mathbf{S}\mathbf{Z}\underline{\bm{\alpha}}^t)]. \label{f1}
\end{align}

In addition, (\ref{4}) and (\ref{5}) can be written as a compact form by collecting all $i \in V$ and $j \in V_i$, which gives $\bm{\omega}^{t+1}=  \underline{\bm{\omega}}^t + \mathbf{S} {\mathbf{Z}} \bm{\alpha}^{t+1}$.
By (\ref{f1}), we have
\begin{align}
\bm{\alpha}^{t+1} =  \arg \min_{\mathbf{n}}  \bigg( & H_r(\mathbf{n}) + \frac{1}{2} \|\mathbf{n} - \bm{\alpha}^t + \mathbf{C} ( \nabla_{\bm{\alpha}} H_s (\bm{\alpha}^t) \nonumber \\
  & + {\mathbf{Z}}^{\top} \underline{\bm{\omega}}^t + \mathbf{Z}^{\top}\mathbf{S}\mathbf{Z}\underline{\bm{\alpha}}^t) \|^2_{\mathbf{C}^{-1}} \bigg),
\end{align}
whose first-order optimality condition can be given by
\begin{align}\label{e-1}
\mathbf{0}  \in & \partial_{\bm{\alpha}} H_r (\bm{\alpha}^{t+1}) + \mathbf{C}^{-1}(\bm{\alpha}^{t+1} - \bm{\alpha}^t)  + \nabla_{\bm{\alpha}} H_s (\bm{\alpha}^t) \nonumber \\
&  + \mathbf{Z}^{\top} \mathbf{S} \mathbf{Z} \underline{\bm{\alpha}}^t + \mathbf{Z}^{\top} \underline{\bm{\omega}}^t  \nonumber \\
= &  \partial_{\bm{\alpha}} H_r (\bm{\alpha}^{t+1}) - \mathbf{C}^{-1}(\bm{\alpha}^t- \bm{\alpha}^{t+1}) + \nabla_{\bm{\alpha}} H_s (\bm{\alpha}^t) \nonumber \\
&  + \mathbf{Z}^{\top} \mathbf{S} \mathbf{Z} \underline{\bm{\alpha}}^t +  \mathbf{Z}^{\top} \bm{\omega}^{t+1}- \mathbf{Z}^{\top} \mathbf{S} \mathbf{Z} \bm{\alpha}^{t+1}.
\end{align}
Noticing that $H_r(\bm{\alpha})$ is convex, we can further have
\begin{align}\label{e2}
H_r  (\bm{\alpha})  -  H_r(\bm{\alpha}^{t+1}) \geq & (\bm{\alpha}-\bm{\alpha}^{t+1})^{\top} \mathbf{C}^{-1} (\bm{\alpha}^t-\bm{\alpha}^{t+1}) \nonumber \\
& - (\bm{\alpha}-\bm{\alpha}^{t+1})^{\top} \nabla_{\bm{\alpha}} H_s(\bm{\alpha}^t)  -(\bm{\alpha}-\bm{\alpha}^{t+1})^{\top} \mathbf{Z}^{\top} \bm{\omega}^{t+1} \nonumber \\
& - (\bm{\alpha}-\bm{\alpha}^{t+1})^{\top} \mathbf{Z}^{\top} \mathbf{S} \mathbf{Z} ( \underline{\bm{\alpha}}^t - \bm{\alpha}^{t+1} ).
\end{align}
Then adding (\ref{e3}) and (\ref{e2}) together gives
\begin{align}\label{r1}
& H (\bm{\alpha}^{t+1}) - H(\bm{\alpha}) \nonumber \\
 \leq & -(\bm{\alpha} -\bm{\alpha}^{t+1})^{\top} \mathbf{C}^{-1} (\bm{\alpha}^t -\bm{\alpha}^{t+1})  \nonumber\\
 & +(\bm{\alpha} -\bm{\alpha}^{t+1})^{\top} \mathbf{Z}^{\top} \bm{\omega}^{t+1} +\| \bm{\alpha}^t - \bm{\alpha}^{t+1}\|^2_{\frac{1}{2}\mathbf{H}} \nonumber\\
&+ (\bm{\alpha} -\bm{\alpha}^{t+1})^{\top} \mathbf{Z}^{\top} \mathbf{S} \mathbf{Z} (\underline{\bm{\alpha}}^t- \bm{\alpha}^{t+1})  \nonumber\\
\stackrel{(i)}{=} & - (\bm{\alpha} -\bm{\alpha}^{t+1})^{\top} \mathbf{C}^{-1} (\bm{\alpha}^t -\bm{\alpha}^{t+1}) \nonumber\\
& - (\bm{\omega} -\bm{\omega}^{t+1})^{\top} \mathbf{S}^{-1} (\underline{\bm{\omega}}^t -\bm{\omega}^{t+1})  \nonumber\\
& - (\bm{\omega}-\bm{\omega}^{t+1})^{\top}\mathbf{Z} \bm{\alpha}^{t+1} +(\bm{\omega}^{t+1})^{\top} \mathbf{Z} \bm{\alpha} \nonumber\\
& -(\bm{\omega}^{t+1})^{\top} \mathbf{Z} \bm{\alpha}^{t+1} + \| \bm{\alpha}^t - \bm{\alpha}^{t+1}\|^2_{\frac{1}{2}\mathbf{H} } \nonumber \\
& + (\bm{\alpha} -\bm{\alpha}^{t+1})^{\top} \mathbf{Z}^{\top} \mathbf{S} \mathbf{Z} (\underline{\bm{\alpha}}^t- \bm{\alpha}^{t+1})  \nonumber\\
\stackrel{(ii)}{=} & \|  \bm{\alpha} -\bm{\alpha}^t \|^2_{\frac{1}{2}\mathbf{C}^{-1}} - \|  \bm{\alpha} -\bm{\alpha}^{t+1} \|^2_{\frac{1}{2}\mathbf{C}^{-1}} \nonumber\\
&  -  \| \bm{\alpha}^t -\bm{\alpha}^{t+1}\|^2_{\frac{1}{2}\mathbf{C}^{-1}}  + \|  \bm{\omega} -\underline{\bm{\omega}}^t \|^2_{\frac{1}{2}\mathbf{S}^{-1} } \nonumber\\
& - \|  \bm{\omega} -\bm{\omega}^{t+1}\|^2_{\frac{1}{2}\mathbf{S}^{-1} } -  \| \underline{\bm{\omega}}^t -\bm{\omega}^{t+1}\|^2_{\frac{1}{2}\mathbf{S}^{-1} }  \nonumber\\
& +(\bm{\omega}^{t+1})^{\top} \mathbf{Z} \bm{\alpha} - \bm{\omega}^{\top} \mathbf{Z} \bm{\alpha}^{t+1} +  \| \bm{\alpha}^t - \bm{\alpha}^{t+1}\|^2_{\frac{1}{2}\mathbf{H}}  \nonumber\\
& - \| \bm{\alpha} -\underline{\bm{\alpha}}^t\|^2_{\frac{1}{2} \mathbf{Z}^{\top}\mathbf{S}\mathbf{Z}} + \| \bm{\alpha} - \bm{\alpha}^{t+1} \|^2_{\frac{1}{2}\mathbf{Z}^{\top}\mathbf{S}\mathbf{Z}} \nonumber \\
& +  \|\underline{\bm{\alpha}}^t -\bm{\alpha}^{t+1}\|^2_{\frac{1}{2}\mathbf{Z}^{\top}\mathbf{S}\mathbf{Z}},
\end{align}
where $(i)$ holds with $ (\bm{\omega} -\bm{\omega}^{t+1})^{\top}( \mathbf{S}^{-1} (\underline{\bm{\omega}}^t - \bm{\omega}^{t+1}) + \mathbf{Z}\bm{\alpha}^{t+1}) = (\bm{\omega} -\bm{\omega}^{t+1})^{\top} \mathbf{S}^{-1} (\underline{\bm{\omega}}^t - \bm{\omega}^{t+1}) +  (\bm{\omega}-\bm{\omega}^{t+1})^{\top} \mathbf{Z}\bm{\alpha}^{t+1} = 0$, and $(ii)$ holds with the fact $\mathbf{m}^{\top} \mathbf{n} = \frac{1}{2} (\| \mathbf{m} \|^2 + \| \mathbf{n} \|^2 - \| \mathbf{m} - \mathbf{n} \|^2) $.

Summing (\ref{r1}) up over $t=0,1,...,T$ from the both sides gives
\begin{align}\label{61}
& \sum_{t=0}^T (H (\bm{\alpha}^{t+1}) - H(\bm{\alpha}) +  \bm{\omega}^{\top} \mathbf{Z} \bm{\alpha}^{t+1}) \nonumber \\
 \stackrel{(i)}{\leq} & \|  \bm{\alpha} -\bm{\alpha}^0 \|^2_{\frac{1}{2}\mathbf{C}^{-1}} - \|  \bm{\alpha} -\bm{\alpha}^{T+1} \|^2_{\frac{1}{2}\mathbf{C}^{-1}} \nonumber\\
& -  \sum_{t=0}^T \| \bm{\alpha}^t - \bm{\alpha}^{t+1} \|^2_{\frac{1}{2}\mathbf{C}^{-1}-\frac{1}{2}\mathbf{H} }   + \|  \bm{\omega} -\bm{\omega}^0\|^2_{\frac{1}{2}\mathbf{Q}}  \nonumber \\
&  - \sum_{t=0}^T \| \underline{\bm{\omega}}^t  -\bm{\omega}^{t+1}\|^2_{\frac{1}{2}\mathbf{S}^{-1}} + \sum_{t=0}^T (\bm{\omega}^{t+1})^{\top} \mathbf{Z} \bm{\alpha}  \nonumber \\
& + \sum_{t=0}^{T} \| \bm{\alpha}^{t+1} - \bm{\alpha}^{t} \|^2_{\frac{1}{2}\mathbf{Z}^{\top}\mathbf{B}\mathbf{Z}} +  \| \bm{\alpha} - \bm{\alpha}^{T+1} \|^2_{\frac{1}{2}\mathbf{Z}^{\top}\mathbf{S}\mathbf{Z}} + \frac{\Gamma}{2}  \nonumber \\
= & \|  \bm{\alpha} -\bm{\alpha}^0 \|^2_{\frac{1}{2}\mathbf{C}^{-1}} - \|  \bm{\alpha} -\bm{\alpha}^{T+1} \|^2_{\frac{1}{2}\mathbf{C}^{-1} - \frac{1}{2}\mathbf{Z}^{\top}\mathbf{S}\mathbf{Z}} \nonumber\\
&  - \sum_{t=0}^T \| \bm{\alpha}^t -\bm{\alpha}^{t+1} \|^2_{\frac{1}{2}\mathbf{C}^{-1}-\frac{1}{2}\mathbf{H} - \frac{1}{2}\mathbf{Z}^{\top}\mathbf{B}\mathbf{Z}} + \|  \bm{\omega} -\bm{\omega}^0\|^2_{\frac{1}{2}\mathbf{Q}} \nonumber \\
& - \sum_{t=0}^T \| \underline{\bm{\omega}}^t  -\bm{\omega}^{t+1}\|^2_{\frac{1}{2}\mathbf{S}^{-1}}  + \sum_{t=0}^T (\bm{\omega}^{t+1})^{\top} \mathbf{Z} \bm{\alpha} + \frac{\Gamma}{2},
\end{align}
where $(i)$ holds by Lemma \ref{ath1}.

Based on the Lagrangian function $\mathcal{L}_H (\bm{\alpha},  \bm{\omega}) = H_s(\bm{\alpha})+ H_r(\bm{\alpha})  + \bm{\omega}^{\top}\mathbf{Z} \bm{\alpha}$, for any $(\bm{\alpha}^*,\bm{\omega}^*) \in {K}$, we have \cite{rockafellar1970convex}
\begin{align}\label{lal2}
\mathcal{L}_H (\bm{\alpha},\bm{\omega}^*) \geq  \mathcal{L}_H (\bm{\alpha}^*,\bm{\omega}^*) \geq  \mathcal{L}_H (\bm{\alpha}^*,\bm{\omega}).
\end{align}
The Karush-Kuhn-Tucker (KKT) conditions are
$\mathbf{0}\in \partial_{\bm{\alpha}} H_r (\bm{\alpha}^*) + \nabla_{\bm{\alpha}} H_s (\bm{\alpha}^*) + \mathbf{Z}^{\top} {\bm{\omega}}^*$ and $\mathbf{Z} \bm{\alpha}^* = \mathbf{0}$ \cite{hanson1981sufficiency}.

Then by letting $\bm{\alpha}=\bm{\alpha}^*$ and $ \bm{\omega} = 2  \frac{\| \bm{\omega}^*\| \mathbf{Z} \overline{\bm{\alpha}}^{T+1}}{\| \mathbf{Z} \overline{\bm{\alpha}}^{T+1} \| }$ in (\ref{61}), we have
\begin{align}\label{61-}
 & (T + 1)(H(\overline{\bm{\alpha}}^{T+1}) - H(\bm{\alpha}^*) +  2\| \bm{\omega}^* \| \|\mathbf{Z} \overline{\bm{\alpha}}^{T+1} \|)  \nonumber\\
\stackrel{(i)}{\leq} & \sum_{t=0}^T \left(H  (\bm{\alpha}^{t+1}) - H(\bm{\alpha}^*) + 2  \frac{\| \bm{\omega}^*\| (\mathbf{Z} \overline{\bm{\alpha}}^{T+1})^{\top}}{\| \mathbf{Z} \overline{\bm{\alpha}}^{T+1} \| } \mathbf{Z} \bm{\alpha}^{t+1} \right) \nonumber\\
\stackrel{(ii)}{\leq}  & \|  \bm{\alpha}^* -\bm{\alpha}^0 \|^2_{\frac{1}{2}\mathbf{C}^{-1}} - \|  \bm{\alpha}^* -\bm{\alpha}^{T+1} \|^2_{\frac{1}{2}\mathbf{C}^{-1} - \frac{1}{2}\mathbf{Z}^{\top}\mathbf{S}\mathbf{Z}} \nonumber\\
&  - \sum_{t=0}^T \| \bm{\alpha}^t -\bm{\alpha}^{t+1} \|^2_{\frac{1}{2}\mathbf{C}^{-1}-\frac{1}{2}\mathbf{H} - \frac{1}{2}\mathbf{Z}^{\top}\mathbf{B}\mathbf{Z}} \nonumber \\
& + \left\|  2  \frac{\| \bm{\omega}^*\| \mathbf{Z} \overline{\bm{\alpha}}^{T+1}}{\| \mathbf{Z} \overline{\bm{\alpha}}^{T+1} \| } -\bm{\omega}^0 \right\|^2_{\frac{1}{2}\mathbf{Q}} \nonumber \\
& - \sum_{t=0}^T \| \underline{\bm{\omega}}^t  -\bm{\omega}^{t+1}\|^2_{\frac{1}{2}\mathbf{S}^{-1}} + \frac{\Gamma}{2}  \nonumber\\
\stackrel{(iii)}{\leq} & \|  \bm{\alpha}^* -\bm{\alpha}^0 \|^2_{\frac{1}{2}\mathbf{C}^{-1}}  + \left\|  2\| \bm{\omega}^*\| \frac{\mathbf{Z} \overline{\bm{\alpha}}^{T+1}}{\| \mathbf{Z} \overline{\bm{\alpha}}^{T+1} \| } -\bm{\omega}^0 \right\|^2_{\frac{1}{2}\mathbf{Q}} + \frac{\Gamma}{2}  \nonumber\\
\stackrel{(iv)}{\leq} & \|  \bm{\omega}^* \|^2_{4\mathbf{Q}} + \| \bm{\omega}^0 \|^2_{\mathbf{Q}} + \|  \bm{\alpha}^* -\bm{\alpha}^0 \|^2_{\frac{1}{2}\mathbf{C}^{-1}} + \frac{\Gamma}{2} \triangleq \Xi,
\end{align}
where $(i)$ is from the convexity of $H$ and $ \sum_{t=0}^T {\bm{\alpha}}^{t+1} = (T+1)\overline{\bm{\alpha}}^{T+1}$, $(ii)$ holds by KKT condition $\mathbf{Z} \bm{\alpha}^*=\mathbf{0}$, $(iii)$ is based on Lemma \ref{lass}, and $(iv)$ holds by {\em{Cauchy-Schwarz inequality}}.

By (\ref{lal2}), we have $H  (\bm{\alpha}^{t+1}) - H(\bm{\alpha}^*) + \bm{\omega}^{*\top} \mathbf{Z} \bm{\alpha}^{t+1} \geq 0$. Letting $\bm{\alpha}=\bm{\alpha}^*$ and $ \bm{\omega} = \bm{\omega}^*$ in (\ref{61}) gives
\begin{align}\label{61+}
 & \sum_{t=0}^T \| \bm{\alpha}^t -\bm{\alpha}^{t+1} \|^2_{\frac{1}{2}\mathbf{C}^{-1}-\frac{1}{2}\mathbf{H} - \frac{1}{2}\mathbf{Z}^{\top}\mathbf{B}\mathbf{Z}} + \sum_{t=0}^T \| \underline{\bm{\omega}}^t  -\bm{\omega}^{t+1}\|^2_{\frac{1}{2}\mathbf{S}^{-1}} \nonumber \\
\leq  & \|  \bm{\alpha}^* -\bm{\alpha}^0 \|^2_{\frac{1}{2}\mathbf{C}^{-1}} + \left\|  \bm{\omega}^* -\bm{\omega}^0 \right\|^2_{\frac{1}{2}\mathbf{Q}} + \frac{\Gamma}{2}.
\end{align}
Then by taking the limit of the both sides of (\ref{61+}) with $T \rightarrow +\infty$, we have $\sum_{t=0}^{+ \infty} \| \bm{\alpha}^t -\bm{\alpha}^{t+1} \|^2_{\frac{1}{2}\mathbf{C}^{-1}-\frac{1}{2}\mathbf{H}- \frac{1}{2}\mathbf{Z}^{\top}\mathbf{B}\mathbf{Z}} < + \infty$ and $\sum_{t=0}^{+ \infty} \| \underline{\bm{\omega}}^t  -\bm{\omega}^{t+1}\|^2_{\frac{1}{2}\mathbf{S}^{-1}} < + \infty$, which means $\lim_{t \rightarrow +\infty} \| \bm{\alpha}^{t+1} - \bm{\alpha}^t \| = 0$ and $\lim_{t \rightarrow +\infty} \| \bm{\omega}^{t+1} - \underline{\bm{\omega}}^t \| = 0$.

By $\bm{\omega}^{t+1}=  \underline{\bm{\omega}}^t + \mathbf{S} {\mathbf{Z}} \bm{\alpha}^{t+1}$, we have $\lim_{t \rightarrow +\infty} \mathbf{Z} \bm{\alpha}^{t+1} = \lim_{t \rightarrow +\infty} \mathbf{Z} \underline{\bm{\alpha}}^{t} = \lim_{t \rightarrow +\infty} \mathbf{S}^{-1}(\bm{\omega}^{t+1}-\underline{\bm{\omega}}^{t}) = \mathbf{0}$ (delay is bounded). Then by taking the limit of the both sides of (\ref{e-1}), we have
\begin{align}\label{}
\mathbf{0}  \in & \lim_{t \rightarrow +\infty} (\partial_{\bm{\alpha}} H_r (\bm{\alpha}^{t+1}) + \mathbf{C}^{-1}(\bm{\alpha}^{t+1} - \bm{\alpha}^t)  + \nabla_{\bm{\alpha}} H_s (\bm{\alpha}^t) \nonumber \\
&  + \mathbf{Z}^{\top} \mathbf{S} \mathbf{Z} \underline{\bm{\alpha}}^t + \mathbf{Z}^{\top} \underline{\bm{\omega}}^t) \nonumber \\
= & \lim_{t \rightarrow +\infty} (\partial_{\bm{\alpha}} H_r (\bm{\alpha}^{t+1}) + \nabla_{\bm{\alpha}} H_s (\bm{\alpha}^t) + \mathbf{Z}^{\top} \underline{\bm{\omega}}^t) \nonumber \\
 & + \lim_{t \rightarrow +\infty}\mathbf{C}^{-1}(\bm{\alpha}^{t+1} - \bm{\alpha}^t) + \lim_{t \rightarrow +\infty}\mathbf{Z}^{\top} \mathbf{S} \mathbf{Z} \underline{\bm{\alpha}}^t  \nonumber \\
= & \lim_{t \rightarrow +\infty} (\partial_{\bm{\alpha}} H_r (\bm{\alpha}^{t+1}) + \nabla_{\bm{\alpha}} H_s (\bm{\alpha}^t) + \mathbf{Z}^{\top} \underline{\bm{\omega}}^t).
\end{align}

To prove the convergence rate, we rearrange (\ref{61-}), which gives
\begin{align}\label{b2}
\frac{\Xi}{T+1}
\stackrel{(i)}{\geq} & H(\overline{\bm{\alpha}}^{T+1})  - H(\bm{\alpha}^*) + 2\| \bm{\omega}^* \| \|\mathbf{Z} \overline{\bm{\alpha}}^{T+1} \| \nonumber \\
{\geq} & H(\overline{\bm{\alpha}}^{T+1})  - H(\bm{\alpha}^*).
\end{align}
Based on (\ref{lal2}) and $\mathbf{Z} \bm{\alpha}^* = \mathbf{0}$, we have
\begin{align}\label{pp1}
& H(\overline{\bm{\alpha}}^{T+1}) - H(\bm{\alpha}^*) \geq - \bm{\omega}^{*\top} \mathbf{Z} \overline{\bm{\alpha}}^{T+1}  \geq - \| \bm{\omega}^* \| \|\mathbf{Z} \overline{\bm{\alpha}}^{T+1} \|.
\end{align}
Adding the both sides of (\ref{b2})-($i$) and (\ref{pp1}) together gives
$\| \bm{\omega}^* \| \|\mathbf{Z} \overline{\bm{\alpha}}^{T+1} \| \leq \frac{\Xi}{T+1}$. Then,
$H(\overline{\bm{\alpha}}^{T+1}) - H(\bm{\alpha}^*) \geq - \| \bm{\omega}^* \| \|\mathbf{Z} \overline{\bm{\alpha}}^{T+1} \|
\geq -\frac{\Xi}{T+1}$ based on (\ref{pp1}). The proof is completed.

\bibliographystyle{IEEEtran}

\bibliography{1myref}
\end{document}